\documentclass[12pt]{amsart}
\pdfoutput=1

\usepackage{cite}

\usepackage{kotex}
\usepackage{lineno,hyperref}
\modulolinenumbers[5]
\usepackage{epstopdf}
\usepackage{pifont}
\usepackage{latexsym}
\usepackage{graphics}
\usepackage{graphicx}
\usepackage{float}
\usepackage[dvips]{xy}
\xyoption{all}
\usepackage{ifpdf}
\usepackage{multirow}
\usepackage{amssymb, amsthm, amsmath, mathrsfs, stmaryrd }
\usepackage{dsfont}
\usepackage{hhline}
\usepackage{centernot}
\usepackage{verbatim, comment, color}
\usepackage{enumerate}
\usepackage{setspace}
\usepackage{lipsum}
\usepackage{subcaption}

\usepackage{mathtools}

\newtheorem{theorem}{Theorem} 
\newtheorem{proposition}[theorem]{Proposition}

\newtheorem{remark}[theorem]{Remark}

\newtheorem{definition}[theorem]{Definition}
\newtheorem{example}[theorem]{Example}

\newcommand{\ba}{\begin{align}}
\newcommand{\ea}{\end{align}}  

\newcommand{\be}{\begin{equation}}

\newcommand{\ee}{\end{equation}}
\newcommand{\bea}{\begin{eqnarray}}
\newcommand{\eea}{\end{eqnarray}}
\newcommand{\barr}{\begin{array}}
\newcommand{\earr}{\end{array}}
\newcommand{\bn}{\begin{enumerate}}
\newcommand{\en}{\end{enumerate}}
\newcommand{\bi}{\begin{itemize}}
\newcommand{\ei}{\end{itemize}}
\newcommand{\bbbm}{\begin{pmatrix}}
\newcommand{\eeem}{\end{pmatrix}}

\newcommand{\cP}{{\cal P}}

\newcommand{\R}{{\mathbb R}}
\newcommand{\N}{{\mathbb N}}
\newcommand{\E}{\mathbb{E}}
\renewcommand{\P}{\mathbb{P}}

\newcommand{\al}{\alpha}

\newcommand{\bt}{\beta}
\newcommand{\ga}{\gamma}

\newcommand{\ignore}[1]{}{}
\newcommand{\no}{\noindent}

\newcommand{\nn}{\nonumber}

\newcommand{\q}{\quad}

\newcommand{{\QED}}{{\hfill QED} \smallskip}

\renewcommand{\subset}{\subseteq}

\renewcommand{\phi}{\varphi}
\newcommand{\cal}{\mathcal}

\DeclareMathOperator{\supp}{supp}
\DeclareMathOperator{\conv}{conv}

 \definecolor{darkspringgreen}{rgb}{0.09, 0.45, 0.27} 
 \definecolor{darkgray}{rgb}{0.66, 0.66, 0.66}
\renewcommand{\sf}{\sffamily}
\numberwithin{equation}{section}
\numberwithin{theorem}{section}
\newcommand{\OT}{{\sf OT }}

\begin{document}
\title
[Vectorial martingale optimal transport]
{Geometry of vectorial martingale optimal transportations and duality
} 

\thanks{The first version of this paper was published on arXiv.org under the title ``{\em Multi-martingale optimal transport}," which was later revised and retitled ``{\em Geometry of multi-marginal martingale optimal transportations and duality}," with the problem referred to as ``{\em Multi-marginal martingale optimal transport problem} (MMOT)". In this version, we refer to the problem as ``{\em Vectorial martingale optimal transport problem} (VMOT)," as suggested by one of the referees, and change the title accordingly. The author gratefully acknowledges the reviewers and the editor for their careful reading, insightful comments, and suggestions.
}

\author{Tongseok Lim}
\address{Tongseok Lim: Mitchell E. Daniels, Jr. School of Business \newline  Purdue University, West Lafayette, Indiana 47907, USA}
\email{lim336@purdue.edu}

\begin{abstract}
The theory of {\em Optimal Transport} (OT) and {\em Martingale Optimal Transport} (MOT) were inspired by problems in economics and finance and have flourished over the past decades, making significant advances in theory and practice. MOT considers the problem of pricing and hedging of a financial instrument, referred to as an option, assuming its payoff depends on a single asset price. In this paper we introduce {\em Vectorial Martingale Optimal Transport} (VMOT) problem, which considers the more general and realistic situation in which the option payoff depends  on multiple asset prices. We address this problem of pricing and hedging given market information -- described by vectorial marginal distributions of underlying asset prices -- which is an intimately relevant setup in the robust financial framework.

We establish that the VMOT problem, as an infinite-dimensional linear programming, admits an optimizer for its dual program. Such existence result of dual optimizers is significant for several reasons: the dual optimizers describe how a person who is liable for an option payoff can formulate optimal hedging portfolios, and more importantly, they can provide  crucial information on the geometry of primal optimizers, i.e. the VMOTs. As an illustration, we show that multiple martingales given marginals must exhibit an extremal conditional correlation structure whenever they jointly optimize the expectation of distance-type cost functions.
 \end{abstract}

\maketitle


\section{Introduction} 
The theory of {\em Optimal Transport} ({\sf OT}) and its probabilistic version {\em Martingale Optimal Transport} ({\sf MOT}) were inspired by problems in economics and finance and have flourished over the past decades, making significant advances in theory and practice.

In this paper we introduce the {\it vectorial Martingale Optimal Transport} ({\sf VMOT}) problem, 
which is an intimately relevant setup in view of the robust financial framework dealing with many asset prices.

To further elaborate this point, let us consider a stochastic process $(X_{t,i})_{t \ge 0}$ of e.g. company $i$'s stock price, where $i \in [d] := \{1,2,...,d\}$. Under the risk neutral probability we may assume $X_t:= \{(X_{t,i})\}_{i \in [d]} \in \R^d$ is a vector-valued martingale. Now in the robust financial framework we shall not assume that the joint probability law of $(X_t)_{t\ge 0}$ is known, since we are unable to determine this joint law from the market information. However, by a standard argument by Breeden and Litzenberger \cite{bl78}, we shall assume that the distribution of each price at each fixed maturity $t_0 >0$, i.e. ${\rm Law}(X_{t_0, i})$, can be observed from the market. 

From now on let us focus on two fixed maturity times $0 < t_1 < t_2$, and denote $X_i := X_{t_1,i}$, $Y_i:= X_{t_2,i}$, $X=(X_1,...,X_d)$, $Y=(Y_1,...,Y_d)$, and assume that $(X,Y)$ is a (one-step, $\R^d$-valued) martingale: $\E[Y|X]=X$. According to the above consideration we do not assume ${\rm Law}(X,Y) \in \cP(\R^{2d})$ is known, but only the $2d$-number of marginal distributions $\mu_i:={\rm Law}(X_i), {\nu_i:=\rm Law}(Y_i)$ are known and fixed. This naturally leads us to consider the following {\sf VMOT} problem: Let $\mu_i, \nu_i \in \cP_1(\R)$, where 
\[
\cP_j(\R^k):= \{ \mu \ | \ \mu \text{ is a probability measure on $\R^k$ with } \int |x|^j d\mu(x) < \infty \}.
\]
Let us write $\vec\mu = (\mu_1,...,\mu_d)$, $\vec\nu= (\nu_1,...,\nu_d)$. We consider the following space of {\em Vectorial Martingale Transportations} from $\vec \mu$ to $\vec \nu$ 
\begin{align}\label{VMT}
{\rm VMT}(\vec\mu,\vec\nu) := \{ \pi \in \cP(\R^{2d}) \ | \  &\pi = {\rm Law} (X,Y), \, \E[Y|X]=X, \\
& X_i \sim \mu_i, \, Y_i \sim \nu_i \q \forall \,  i \in [d] \}. \nn
\end{align}
Let $c : \R^{2d} \to \R$ be a (cost) function. We define the {\sf VMOT} problem as
\begin{align}\label{VMOT}
\max / \, {\rm minimize } \ \ \E_\pi [c(X,Y)] \ \text{ over } \ \pi \in {\rm VMT}(\vec\mu, \vec\nu).
\end{align}
The function $c$ has a natural interpretation in finance as an exotic option whose payoff is determined by the $d$-number of asset prices $(X,Y)$ at the terminal maturity. $ \E_\pi [c(X,Y)]$ shall be considered as a fair price of the option $c$ in the risk-neutral world governed by the law of asset prices $\pi$, but since $\pi$ would not be observed from the financial market, we need to consider all possible laws ${\rm VMT}(\vec\mu,\vec \nu)$ given the marginal information $\vec\mu, \vec\nu$. Under this information, the max / min value in \eqref{VMOT} can be interpreted as the upper / lower price bound for $c$ respectively.

Note that without the martingale constraint $\E[Y|X]=X$, \eqref{VMOT} would be the usual multi-marginal optimal transport problem. But \eqref{VMOT} still belongs to the class of {\em infinite-dimensional linear programming}, since the martingale constraint is linear. Now unlike the optimal transport case, the problem \eqref{VMOT} will be feasible (i.e. ${\rm VMT}(\vec\mu, \vec\nu) \neq \emptyset$) if only if every pair of marginals $\mu_i, \nu_i$ is in convex order, defined by
\[ \mu_i \le_c \nu_i \ \ \text{iff} \ \ \mu_i, \nu_i \in \cP_1 \ \text{ and } \ \int f d\mu_i \le \int f d\nu_i \ \ \text{for every convex } f,
\]
as shown by Strassen \cite{St65}. Thus we will always assume $\mu_i \le_c \nu_i$ for all $i \in [d]$ in {\sf VMOT} problem. And, because we will need to calculate means and potential functions of the measures, we will make the following assumption throughout the paper (unless otherwise specified):
\vspace{1mm}

\noindent {\bf Assumption.} All measures appearing in this paper 
are assumed to have finite first moments.
\vspace{1mm}

Our first main result investigates the extremal correlation behavior of the terminal prices $Y$ conditional on $X$, when the martingales $(X,Y)$ jointly max/minimize a distance-type cost function given vectorial marginals. More precisely, we establish the following theorem.
\begin{theorem} \label{main1} Assume  all $\mu_i$ are absolutely continuous with respect to Lebesgue measure, denoted as $\mu_i \ll \cal L$. If $\pi$ solves the maximization problem in \eqref{VMOT} with the Euclidean distance cost $c(x,y)= |x-y|$, then given $X$, the conditional distribution of $Y$ lies in the set of extreme points of a convex set $C_X \subset \R^d$ depending on $X$, $\pi$--almost surely.
\end{theorem}
Here we say that a random variable $Y$ lies in $A$ if $Y \in A$ a.s., and a distribution $\xi$ lies in $A$ if $\xi$ has its full mass in $A$. We shall also prove an analogous result for the minimization problem in \eqref{VMOT}, and we also note that Theorem \ref{main1} still holds for more general class of strictly convex costs $c(x,y)=||x-y||$. See Theorem \ref{extreme} for more details. We further note that the statement of extremal correlation was firstly given and studied by \cite{GKL2} in the context of {\sf MOT} and two $n$-dimensional marginals.
\\

Recall that the {\sf VMOT} problem belongs to the class of infinite-dimensional linear programming, which implies that the problem has its {\em dual programming} problem. For the (primal) minimization problem in \eqref{VMOT}, its dual problem is given by 
 \begin{equation}\label{dualproblem}
 \sup_{(\phi_i,\psi_i,h_i) \in \Psi}\sum_{i=1}^d  \big{( } \int \phi_{i}d\mu_i + \int \psi_{i}d\nu_i \big{)}
 \end{equation}
where $\Psi$ consists of triplets $\phi_i, \psi_i : \R \to \R \cup \{-\infty\}$ and $h_i : \R^d \to \R$ such that $\phi_i \in L^1(\mu_i)$, $\psi_i \in L^1(\nu_i)$, $h_i$  is bounded for every $i \in [d]$, and 
\begin{align}
\label{ptwiseineq} &\sum_{i=1}^d \big{(} \phi_{i} (x_i) + \psi_{i} (y_i) + h_{i}(x) (y_i - x_i) \big{)} 
 \le c(x,y) \q \forall (x,y) \in \R^{2d}
\end{align}
where $x=(x_1,.,,,x_d)$, $y=(y_1,...,y_d)$. Analogous dual problem for the maximization in \eqref{VMOT} is easily obtained by changing sign of $c$ to $-c$. 

The dual problem also has an important interpretation in finance. Suppose a financial firm is liable for an option $c$ it has sold so that the firm has to pay $c(X,Y)$ at the terminal maturity $t_2$. To hedge its risk, the company may consider buying European options $\phi_i, \psi_j$ where each payoff relies only on a single asset at a fixed time. In addition, the firm may consider holding $h_{i}$ number of shares of the $i^\textrm{th}$ asset held between the $t_1$ and $t_2$ maturity, such that its payoff at $t_2$ is $h_i(X)\cdot(Y_i-X_i)$. Notice $h_{i}$ is a function of the past prices of all assets $\{X_i\}_{i=1}^d$. Thus the left hand side of \eqref{ptwiseineq} yields the overall payoff of the hedging portfolio $(\phi_i, \psi_i, h_i)_{i=1}^d$, and the inequality \eqref{ptwiseineq} imposes that the position must subhedge the liability for all scenarios of $(X,Y)$. Now notice the dual of the maximization in \eqref{VMOT} becomes an optimal superhedging problem.\\

Now the celebrated {\em duality} result asserts, under a mild assumption on $c$, that the primal and dual optimal values coincide (see e.g. \cite{eglo, Z})
\begin{align}\label{duality}
P(c):&=\inf_{\pi \in {\rm VMT}(\mu, \nu)} \E_\pi [c(X,Y)] \\
&= \sup_{(\phi_i,\psi_i,h_i) \in \Psi}\sum_{i=1}^d  \big{( } \int \phi_{i}d\mu_i + \int \psi_{i}d\nu_i \big{)} =:D(c). \nn
\end{align}

Still a natural question is whether a dual optimizer exists, that is, a ``tightly" sub/superhedging portfolio can be constructed. In the past decades, researchers working in {\sf OT} already realized that the existence of an optimizer to the dual problem -- often called {\em dual attainment} -- is  the essential key for the investigation of the optimal transport plans. But due to the nature of infinite dimensionality, establishing the dual attainment turns out to be quite subtle, as can be seen e.g. by the work of Brenier \cite{br}. Thus it is no surprise that our second main result on the dual attainment is crucial for establishing Theorem \ref{main1}. For a measure $\xi$ on $\R$ define its {\em potential function} $u_\xi (x) := \int |x-y| d\xi(y)$. We say 
$\mu\leq_{c}\nu$ is \emph{irreducible} if $I:=\{u_{\mu}< u_{\nu}\}$ is connected and $\mu(I)=\mu(\R)$. For more information on the notion of irreducibility, see Section \ref{m2} or \cite{bnt}.
\begin{theorem} \label{main2}
Let $(\mu_i,\nu_i)_{i \in [d]}$ be irreducible pairs of marginals on $\R$. Let $c(x_1,...,x_d,y_1,...,y_d)$ be a lower-semicontinuous cost such that $|c(x,y)| \le \sum_{i=1}^d \big{(}v_i(x_i) + w_i(y_i)\big{)}$ for some continuous $v_i \in L^1(\mu_i)$, $w_i \in L^1(\nu_i)$. 
Then there exists a {\em pointwise dual maximizer} {\rm({\sf PDM})},  a triplet of functions $(\phi_i, \psi_i, h_i)_{i=1}^d$ that satisfies \eqref{ptwiseineq} tightly in the following pointwise manner (but does not necessarily belong to $\Psi$):
\begin{align}
 \label{ptwiseeq}& \sum_{i=1}^d \big{(} \phi_{i} (x_i) + \psi_{i} (y_i) + h_{i}(x) (y_i - x_i) \big{)} 
 = c(x,y) \q \pi-a.s.
\end{align}
for every {\rm VMOT} $\pi$ which solves the minimization problem in \eqref{VMOT}.
\end{theorem}
 We analogously define a {\em pointwise (or pathwise) dual minimizer} ({\sf PDm}) for which the inequality in \eqref{ptwiseineq} should be reversed and $\pi$ in \eqref{ptwiseeq} now solves the maximization problem in \eqref{VMOT}. This simply corresponds to the change of $c$ to $-c$ in the theorem. We shall sometimes call Theorem \ref{main2} as {\em pointwise/pathwise dual attainment}  ({\sf PDA}) for {\sf VMOT} problems.

The term {\em pointwise} indicates the fact that \eqref{ptwiseineq}, \eqref{ptwiseeq} hold in a pointwise manner (i.e., the equality \eqref{ptwiseeq} is satisfied for $\pi$ - almost every ``point" or ``path" $(x,y)$), and that we do {\em not} impose the integrability condition, namely $\phi_i \in L^1(\mu_i)$, $\psi_i \in L^1(\nu_i)$, or $h_i$  is bounded, on the {\sf PDM}. Thus, a  {\sf PDM} needs not be in $\Psi$, and moreover, $\phi_i, \psi_i$ may assume the value $-\infty$. However, it is indeed true that $\phi_i \in \R$ $\mu_i$-a.s. and $\psi_i \in \R$ $\nu_i$-a.s., implied by the equation \eqref{ptwiseeq}. See also Remark \ref{remark3.2}.

A couple of pioneering works on the dual attainment for {\sf MOT} problems on the line ($d=1$) was established by Beiglb{\"o}ck-Juillet \cite{bj} and Beiglb{\"o}ck-Nutz-Touzi \cite{bnt}, where they already observed that the {\sf PDA} could fail if one insists the integrability condition on {\sf PDM/PDm}, or if one does not impose the irreducibility on the marginals. Beiglb{\"o}ck-Lim-Ob{\l}{\'o}j \cite{blo} established {\sf PDA} without the irreducibility assumption but imposing further regularity condition on $c$. Theorem \ref{main2} is a continuation of these efforts for the case $d \ge 2$, i.e., for the {\sf VMOT} setup.\\

\no{\bf Further background.} Let us briefly introduce the {\it optimal transport} ({\sf OT}) problem, which is the prototype of {\sf MOT} and {\sf VMOT} problems. To begin with, let us introduce the notion of {\em push-forward} of measures. Given a measurable map $F : \cal X \to \cal Y$ and a measure $\mu$ on $\cal X$, the push-forward of $\mu$ by $F$, denoted by $F_\#\mu$, is a measure on $\cal Y$ satisfying
\[ 
F_\#\mu(A) = \mu(F^{-1}(A))
\]
for every measurable $A \subset \cal Y$.

Let $p^x, p^y : \R^d \times \R^d \to \R^d$ be the projection maps given by $p^x(x,y)=x$, $p^y(x,y)=y$. For a given cost function $c:\R^d\times \R^d \to \R$ and two Borel probability measures $\mu, \nu$ on $\R^d$, the {\sf OT} problem is given by
\begin{equation}\label{TP}
\mbox{Minimize} \quad \text{cost}[\pi] = \int_{\R^d\times \R^d} c(x,y) \,d\pi(x,y) \quad\mbox{over}\quad \pi\in \Pi(\mu,\nu)
\end{equation}
where $\Pi(\mu,\nu) \subset \cP(\R^{2d})$ is the set of {\it  transport plans}, or {\it couplings}, if all $\pi \in \Pi(\mu,\nu)$ have given marginals $\mu$ and $\nu$, meaning $p^x_\#\pi = \mu$, $p^y_\#\pi = \nu$.

Kantorovich \cite{K1,K2} proposed the general formulation of the problem \eqref{TP} in the 1940s after Gaspard Monge's first formulation in the 1780s \cite{Mo1781}. Since then, many important contributions have been made on the subject, as can be seen by part through the works \cite{akp,ap,bc,br,cfm, cms01, cms06, evans,gm, MTW, Mc97, Sa15,TW2,Vi03,Vi09}. In particular, one of the central questions is when the {\sf OT} plans are given by transport maps -- often called {\em Monge problem}, especially when $c(x,y)=|x-y|$ -- that is, when there exists a map $T : \R^d \to \R^d$ such that a minimizer $\pi$ in \eqref{TP} is given by $\pi = ({\rm Id}, T)_\#\mu$, where ${\rm Id}(x)=x$ is the identity map on $\R^d$.

To answer, we need to understand the geometry of the {\sf OT} plans. Such a characteristic structure of {\sf OT} plan  is encoded in its support, which exhibits the so called {\it $c$-cyclical monotonicity}. The monotonicity is already useful, but when combined with the celebrated theorem of Rockafellar \cite{Rock}, we obtain the {\em dual attainment}, a cornerstone of {\sf OT}: there exist two functions $\phi, \psi :\R^d \to \R$, called {\em Kantorovich potentials}, such that for every minimizer $\pi$ of the problem \eqref{TP}, we have
\begin{align}
 \label{otineq} &\phi(x)+ \psi(y)  \le c(x,y)  \quad  \forall x \in \R^d,\, \forall y \in \R^d,\\
 \label{oteq} &\phi(x)+ \psi(y)  = c(x,y)  \quad \,\, \,\pi - a.e.\, (x,y).
\end{align}
It turns out the dual attainment is  essential for understanding {\sf OT}\cite{gm}.

Let us view the \OT  maps in a slightly different angle, and for this we recall the notion of {\em disintegration} of measures. For a general measure $\pi \in \cP(\R^{2d})$, a disintegration of $\pi$ via its first marginal $\pi^1:=p^x_\# \pi$ is 
\be
\pi = \pi_x \otimes \pi^1
\ee
which means for any Borel sets $A, B \subset \R^d$, $\pi(A \times B) = \int_A \pi_x(B) d\pi^1(x)$. The family of probabilities $\pi_x := \{\pi_x\}_x$ is called a disintegration (or  {\em kernel}) of $\pi$ with respect to $\pi^1$. It is known this family uniquely exists $\pi^1$- a.e. $x$. Perhaps probabilistic language helps to understand better. Suppose $\pi$ is the joint distribution of the $\R^d$-valued random variables $X, Y$, denoted as $\pi = {\rm Law}(X,Y)$. Then $\pi^1 = {\rm Law}(X)$, and $\pi_x$ is the conditional law of $Y$ given $X=x$, that is  $\pi_x(B) = \P(Y \in B \, | \, X=x )$.

Now ``$\pi \in \Pi(\mu,\nu)$ is given by a transport map $T$." is equivalent to
\be
\pi = ({\rm Id}, T)_\#\mu \iff \pi = \pi_x \otimes \mu \ \text{ with } \ \pi_x = \delta_{T(x)}
\ee
where $\delta_y$ is the {\em Dirac mass} at $y$, the probability measure concentrated at $y$. In this case one could say  the solution $\pi$ exhibits an extremal geometry in the sense that a disintegration $\pi_x$, which describe the solution, are as extreme as Dirac masses, the most condensed measures.
Finally, we would like to also mention some important applications of {\sf OT} in economics \cite{fkm11, dlms13, p12, cghp21}, and in statistical theory \cite{ccg16, cghh17, ru, cfg10}.
\\

Recently a variant of {\sf OT}, referred to as  {\em martingale optimal transport} ({\sf MOT}), was introduced. In {\sf MOT} problem,  we consider the following:
\begin{equation}\label{MTP}
\mbox{Minimize}\ \ \text{Cost}[\pi] = \int_{\R^d \times \R^d} c(x,y)\,d\pi(x,y)\ \ \mbox{over}\ \  \pi\in {\rm MT}(\mu,\nu)
\end{equation}
where MT$(\mu,\nu)$ (Martingale Transport plans) is a subset of $\Pi(\mu,\nu)$ such that each $\pi = \pi_x \otimes \mu \in$ MT$(\mu,\nu)$ satisfies $\delta_x \le_c \pi_x$, that is $\pi_x$ has its {\it barycenter} at $x$. Notice the kernel $(\pi_x)_x$ cannot be Dirac masses unless $\pi_x = \delta_x$. Nonetheless, Theorem \ref{main1} shows the kernel of each {\sf MOT} for $l_1$ norm-type costs still exhibits an interesting extremal distribution.

Probabilistic description of the {\sf MOT} problem is as follows: consider
\begin{align}\label{opt}
\text{Minimize}\quad  
\E_{\P} \,[c(X,Y)]
\end{align}
over all couples $(X,Y)$ of $\R^d$-valued random variables on some probability space $(\Omega, {\mathcal F}, \P)$ such that ${\rm Law}(X)=\mu$, ${\rm Law}(Y)=\nu$ and $\E[Y|X]=X$, i.e. $(X,Y)$ is a one-step martingale. Strassen \cite{St65} showed  MT$(\mu, \nu)$ is nonempty if and only if $\mu$, $\nu$ are in convex order $\mu \le_c \nu$.\\

 While some pioneering works to investigate the {\sf MOT} problem have been made including \cite{BeHePe11,ds1,GaHeTo11,HoKl12,HoNe11}, there is a related problem, called the {\it Skorokhod embedding problem} ({\sf SEP}) which has a long history in probability theory.  Since Hobson \cite{Ho98,Ho11} recognised the important connection between {\em model independent finance and asset pricing theory} and {\sf SEP} (see Ob{\l}{\'o}j \cite{Obloj} for a nice overview of the {\sf SEP} and Beiglb{\"o}ck-Cox-Huesmann \cite{bch} for a link with {\sf OT} theory), much related research  has been done in this context including \cite{cko, gtt1, gtt2, os, cot19} for instance.
 
In the early stages of development of {\sf MOT} theory, most of the research focused on single asset cases. More recently there have been efforts to generalize the theory into higher dimension $d \ge 2$, especially around the theory of duality and its attainment \cite{Lim, GKL2,GKL3, d18, dt19, eglo, os17}. Nonetheless, to the best of the author's knowledge a complete understanding of the possibility for the dual attainment has not been achieved: the added ``trading strategy" term $h(x) \cdot (y-x)$  in the dual formulation seems to drastically change the picture from the {\sf OT} case, and the potential non-existence of dual optimizers makes the study of solutions to {\sf MOT} much more complex in the multi-dimensional setting.
\\

While the {\sf MOT} formulation \eqref{MTP} is a natural analogue of {\sf OT}, not only it has difficulty in establishing {\sf PDA}, but also its assumption on the marginal information is rather unrealistic from the financial point of view. This is because {\sf MOT} assumes the {\em joint distribution} of the assets $\{X_i\}_i$ and $\{Y_i\}_i$ are known -- the $\mu,\nu$. However, in practice such information is hardly observable from the market, and this is part of the motivation for us to introduce the current {\sf VMOT} formulation. 
Let us recall that the main difference between {\sf MOT} and {\sf VMOT}  lies in the marginal constraint, since in {\sf VMOT} only the individual laws of $\{X_i\}_i$ and $\{Y_i\}_i$ -- the $\vec \mu, \vec \nu$ -- are assumed to be known. Thus, there are now at least three unknown couplings in the {\sf VMOT} problem whose structures are of major interest:  the {\sf VMOT} plan $\pi = {\rm Law}(X,Y)$, and moreover the induced couplings $\pi^1:=p^x_\#\pi =  {\rm Law}(X)$ and $\pi^2 := p^y_\#\pi=  {\rm Law}(Y)$.
\\

This paper is organized as follows. In section \ref{m1} we state and prove a more detailed version of Theorem \ref{main1}  by applying {\sf PDA}. In section \ref{m2} we prove the {\sf PDA}, Theorem \ref{main2}, using a compactness result whose proof is given in an appendix. Section \ref{m3} presents further results on the structure of {\sf VMOT} and optimality of induced couplings.

\section{Extremal correlation of martingales -- Theorem \ref{main1}} \label{m1}

In this section, we state and prove Theorem \ref{extreme} which is a more precise and extended version of Theorem \ref{main1}. We will prove it here by assuming Theorem \ref{main2}. Let us introduce some terminology. A norm $|| \cdot||$ on $\R^d$ is called {\em strictly convex} 
if its unit ball $B=\{x : ||x|| \le 1 \}$ is strictly convex, i.e. every boundary point of $B$ is an extreme point of $B$. Euclidean norm $| \cdot |$ clearly belongs to this class. For a set $A$, $\conv(A)$ is the convex hull of $A$, and for a convex set $A$, ${\rm Ext}(A)$ is the set of extreme points of $A$. For a measure $\xi$, the {\em support} of $\xi$, denoted by $\supp \xi$, is the smallest closed set on which $\xi$ has its full mass. For measures $\mu,\nu$, there exists a unique largest measure $\mu \wedge \nu$ (called the {\em common mass} of $\mu,\nu$) which is dominated by $\mu$ and $\nu$, i.e. $\mu \wedge \nu \le \mu$, $\mu \wedge \nu \le \nu$, where $\mu \le \nu$ means $\mu(A) \le \nu(A)$ for every measurable $A$. When $\mu,\nu$ are given by densities $f,g$ respectively, $\mu \wedge \nu$ is given by the density $f \wedge g := {\rm min}(f,g)$. $\mu \ll \cal L^d$ means $\mu$ is absolutely continuous w.r.t. the Lebesgue measure on $\R^d$, that is $\mu$ has a density function.

\begin{theorem} \label{extreme} Let $(\mu_i,\nu_i)_{i \in [d]}$ be pairs of probability measures on $\R$ in convex order, and $\mu_i \ll \cal L^1$ for all $i$. Let $c(x,y) = \pm || x-y||$ where the norm $|| \cdot ||$ is strictly convex 
and  $x \mapsto ||x||$ is differentiable on $ \R^d \setminus \{0\}$. Let  $\pi = \pi_x \otimes \pi^1$ be any {\em minimizer} for  \eqref{VMOT} where $\pi^1= p^x_\# \pi$, $\pi^2= p^y_\# \pi$.
\begin{enumerate}
\item If $c(x,y) = - || x-y||$, then the support of $\pi_x$ coincides with the extreme points of the  convex hull of itself:
\[
{\rm supp}\, \pi_x = {\rm Ext}\,\big{(} {{\rm conv}}({\rm supp}\, \pi_x) \big{)}, \quad \pi^1 - a.e. \, x.
\]
\item If $c(x,y) =  || x-y||$, we have $D_\#( \pi^1 \wedge \pi^2) \le \pi$ where $D(x)=(x,x)$, meaning that the common marginal $\pi^1 \wedge \pi^2$ stays put under $\pi$, or the minimizer $\pi$ does not move the mass of $\pi^1 \wedge \pi^2$. Let $\tilde \pi := \pi - D_\# (\pi^1 \wedge \pi^2)$ and disintegrate as $\tilde \pi = \tilde \pi_x \otimes \tilde \pi^1$. Then 
\[
{\rm supp}\, \tilde \pi_x = {\rm Ext}\,\big{(} {{\rm conv}}({\rm supp}\, \tilde\pi_x) \big{)}, \quad \tilde\pi^1 - a.e.\, x.
\]
\end{enumerate}
\end{theorem}
Observe $\tilde \pi$ is a martingale transport from $\tilde \pi^1:=p^x_\# \tilde \pi = \pi^1 - \pi^1 \wedge \pi^2$ to $\tilde \pi^2 :=p^y_\# \tilde \pi = \pi^2 - \pi^1 \wedge \pi^2$, and $\tilde \pi^1 \wedge \tilde \pi^2 = 0$. Since $D_\#( \pi^1 \wedge \pi^2)$ is an identity transport, we conclude that by minimizing $\E||x-y||$, we get
\[
{\rm supp}\, \pi_x = {\rm Ext}\,\big{(} {{\rm conv}}({\rm supp}\, \pi_x) \big{)} \ \text{ or } \ {\rm supp}\, \pi_x = {\rm Ext}\,\big{(} {{\rm conv}}({\rm supp}\, \pi_x) \big{)} \cup \{x\}
\]
and moreover the latter will be the case for $\pi^1 \wedge \pi^2$ -- a.e. $x$.

We remark that \cite{HoKl12}, \cite{HoNe11} studied the structure of {\sf MOT} w.r.t. $c(x,y) = \pm |x-y|$ in the single-martingale setup, and later \cite{bj} rediscovered part of their results from a different approach. Theorem \ref{extreme} may be viewed as a generalization of e.g. \cite[Theorem 7.3, Theorem 7.4]{bj} to the multiple martingales setup, showing that the martingales given marginals are correlated in an extremal way to achieve the optimum.

\begin{proof}
{\bf Step 1.} 
Let us tentatively assume that $(\mu_i,\nu_i)$ are  irreducible for every $i$. This assumption will be removed in the last step. Then by Theorem \ref{main2}  we have a {\sf PDM} $(f_i, -g_i, h_i)_{i \in [d]}$, that is, by denoting $ f^\oplus(x) = \sum_{i=1}^d f_i(x_i)$, $g^\oplus(y) = \sum_{i=1}^d g_i(y_i)$ and $h=(h_1,...,h_d)$, we have

\begin{align}
\label{ptwiseineq1} f^\oplus(x)+ h(x)\cdot (y-x) &\le c(x,y) +g^\oplus(y) \q \forall (x,y) \in \R^{2d}, \text{ and} \\
 \label{ptwiseeq1} f^\oplus(x)+ h(x)\cdot (y-x) &= c(x,y) +g^\oplus(y) \q \pi-a.s..
\end{align}
Recall that  in Theorem \ref{main2}, \eqref{ptwiseeq1} is attained as real-valued, yielding $f_i$ is finite $\mu_i$-a.s. and $g_i$ finite $\nu_i$-a.s.. In this step we shall obtain a differential identity \eqref{nablazero}.

To this end, recall the martingale Legendre transform of $g^\oplus(y)$ is a pair of functions $(\alpha,\tilde \gamma)$ defined as in \cite[Definition 3.1]{GKL2}:
\begin{align}\label{alpha1}
 \alpha (x) := \sup \{&a \in \R \,|\, \text{there exists } b \in \R^d \, \text{ such that } \\
&a + b \cdot (y-x)  \le c(x,y)+g^\oplus(y) \text{\, for all } y\},  \nn\\
\label{gamma1}
\tilde \gamma(x) := \{b \in \R^d \,|\, &\alpha(x) + b \cdot (y-x) \le  c(x,  y)+g^\oplus(y) \text{\, for all } y \}.  
\end{align}
In general $\tilde \gamma$ is a convex set-valued (possibly empty, but see \eqref{regularity1} below) function and we may choose a $\R^d$-valued function $\gamma \in \tilde \gamma$. Observe that if we define $H_x(y)$ to be the convex envelope of $y \mapsto c(x,y) + g^\oplus(y)$, then $y \mapsto \alpha(x) +\gamma(x) \cdot (y-x)$ is an affine tangent function to $H_x(y)$ at $x$. Now recall that $\mu_i\leq_{c}\nu_i$ is called \emph{irreducible} if $I_i:=\{u_{\mu_i}< u_{\nu_i}\}$ is connected open and $\mu_i(I_i)=\mu_i(\R)$. Set $I := \otimes_i I_i \subset \R^d$ to be an open rectangle. Now \eqref{alpha1} gives $ f^\oplus(x) \le \alpha(x)$, thus
 \begin{align*}
 f^\oplus(x)+ \gamma(x)\cdot (y-x) - c(x,y)& \le \alpha(x) + \gamma(x)\cdot (y-x) - c(x,y) \\
 &\le g^\oplus(y)  
\end{align*}
and equality holds throughout (as real-valued) $\pi$ - a.s.. The fact $g_i$ finite $\nu_i$-a.s. then yields $I \subset \conv\{ g^\oplus \text{ is finite}\}$. With this and the fact that $c$ is Lipschitz, \cite[Theorem 3.2]{GKL2} yields the following local regularity:
\begin{align}\label{regularity1}
\alpha \text{ is Lipschitz and }  \gamma \text{ is bounded on every compact subset of }  I.
\end{align}
We wish $f_i$ to attain Lipschitz property, so we will take the Legendre transform of $f_i$ with respect to $\alpha$ to obtain such property. But as $\alpha$ is in general only locally Lipschitz, let us take the transform locally as follows:  write $I_i = ]a_i, b_i[$ and for $n \in \N$, let $I^n_i =  ]a_i + \frac1n, b_i - \frac1n[ \cap ]-n,n[$ and $I^n = \otimes_i I^n_i$. Now define $\phi_i$ successively for $i=1,2,...,d$ by
\begin{align}\label{partiallegendre}
\phi_i(x_i) := \inf_{x_j \in I^n_j, j \neq i} \big{(} \alpha(x) - \sum_{j <  i} \phi_j(x_j)  - \sum_{j >  i} f_j(x_j) \big{)}.
\end{align}
Since $\alpha$ is (globally) Lipschitz in $I^n$ (because $\al$ is locally Lipschitz in $I$), \eqref{partiallegendre} implies that
$\phi_i$  is   Lipschitz  in $ I_i^n$, and 
\begin{align*}
f_i \le \phi_i  \, \text{ on } \, I^n_i \quad \text{and} \quad \sum_i \phi_i(x_i) \le \alpha(x)  \, \text{ on } \, I^n.
\end{align*}

Let $\pi_n := \pi  \big{|}_{ I^n \times \R^d}$ be the restriction of $\pi$ on $I^n \times \R^d$, let $\pi_n^1=p^x_\#\pi_n$, and let $(\mu^n_i)_{i\in[d]}$ be the one-dimensional marginals of $\pi_n^1$.  As $\mu^n_i \le \mu_i $ yields $\mu^n_i \ll \cal L$, each $\varphi_i$ is differentiable $\mu^n_i$ -a.e.. Hence 
\begin{align}\label{phidiff}
\phi^\oplus(x) := \sum_i \phi_i(x_i) \,\text{ is differentiable} \quad \pi^1_n - a.e. \, x.
\end{align}
Now because $f^\oplus \le \phi^\oplus \le \al$, we have
\begin{align}
 \label{epsilonineq} &\phi^\oplus(x)+ \gamma(x)\cdot (y-x) - c(x,y)  \le g^\oplus(y)  \quad  \forall (x,y) \in I^n \times \R^d,\\
 \label{epsiloneq} &\phi^\oplus(x)+ \gamma(x)\cdot (y-x) - c(x,y)  = g^\oplus(y)  \quad  \pi_n - a.e.
\end{align}
We may rewrite \eqref{epsiloneq} as
\begin{align}\label{xeq}
\text{For } \pi_n^1 - a.e. \,x, \ \phi^\oplus(x)+ \gamma(x)\cdot (y-x) - c(x,y)  = g^\oplus(y)  \ \  \pi_x - a.e.\,y.
\end{align}
Fix $x_0$ at which \eqref{xeq} holds and $\phi^\oplus$ is differentiable. Let $V_0$ be the subspace of $\R^d$ spanned by the set $\supp \pi_{x_0}  - x_0$ (translation of $\supp \pi_{x_0}$ by $-x_0$). If dim $V_0 = 0$ then it simply means  $\supp \pi_{x_0}  = \{x_0\}$ and there is nothing to prove. Thus let us assume that dim $V_0 \ge 1$. Now \cite[Lemma 4.1]{GKL2} showed that, if \eqref{epsilonineq}, \eqref{xeq} hold, then
\begin{align}\label{gammadiff}
\text{directional derivative of } {\rm proj}_{V_0} \gamma \,\text{ exists at $x_0$ in every direction $u \in V_0$, }
\end{align}
where ${\rm proj}_{V_0} \gamma$ is the orthogonal projection of the $\R^d$-valued function $\gamma$ on  $V_0$. While a proof can be found in \cite{GKL2}, here let us give some intuition: in \eqref{epsilonineq} the function $y \mapsto \phi^\oplus(x_0)+ \gamma(x_0)\cdot (y-x_0) - c(x_0,y)$ is bounded above by $g^\oplus(y)$, but \eqref{xeq} tells us that the function is in fact tightly bounded by $g^\oplus(y)$ for $\pi_{x_0}$ - a.e.\,$y$. The facts $\pi_{x_0}$ has barycenter at $x_0$ and the functions $\varphi^\oplus(x)$ and $x \mapsto c(x,y)$ are differentiable at $x = x_0$ implies that ${\rm proj}_{V_0} \gamma$ has not much room to vary wildly near $x_0$ in $V_0$, yielding it is forced to be differentiable at $x_0$ in $V_0$.

Next, note that \eqref{epsilonineq}, \eqref{xeq}  imply for $\pi_{x_0} - a.e.\,y$
\begin{align}\label{goodineq}
&\phi^\oplus(x_0)+ \gamma(x_0)\cdot (y-x_0) - c(x_0,y) \\
&\ge \phi^\oplus(x)+ \gamma(x)\cdot (y-x) - c(x,y) \quad \forall x \in I^n, \nn
\end{align}
and notice by continuity of $c$, \eqref{goodineq} holds for every $y \in \supp \pi_{x_0}$. Then from \eqref{phidiff},  \eqref{gammadiff}, we deduce that for any nonzero vector $u$ in $V_0$, by taking $u$-directional derivative in $x$ at $x_0$, it holds that for any $y \in \supp \pi_{x_0} \setminus \{x_0\}$,
\begin{align}\label{nablazero}
\nabla_u \phi^\oplus(x_0)+ \nabla_u \gamma(x_0) \cdot (y-x_0) - \gamma(x_0) \cdot u  - \nabla_u c(x_0,y) = 0. 
\end{align}
From this identity \cite{GKL2} proved the theorem when the cost is given by the Euclidean norm. We will follow  a similar  line but as we deal with more general strictly convex norms, we shall need a more involved argument.\\

{\bf Step 2.}  We shall prove the ``non-staying" property of the common mass $\pi^1 \wedge \pi^2$ in the case $c(x,y) = -||x-y||$ ($\pi^1 \neq \pi^2$ by irreducibility) and the ``staying" property in the case $c(x,y) = +||x-y||$ for the minimization problem in \eqref{VMOT}.

For $c(x,y) = - ||x-y||$,  \eqref{goodineq} immediately implies $x_0 \notin \supp \pi_{x_0}$ for every $x_0$ at which $\phi^\oplus$ is differentiable, thus $\pi^1$-a.s.,
 (which we call the non-staying property of the minimizer $\pi$) by the following reason: if $x_0 \in \supp \pi_{x_0}$, then the function $x \mapsto  \phi^\oplus(x)+ \gamma(x)\cdot (x_0-x) + ||x - x_0||$ must attain its maximum at $x=x_0$ by  \eqref{goodineq}. But notice that due to the increase of $x \mapsto ||x-x_0||$ the function will strictly increase as $x$ moves away from $x_0$ along any direction $u \in V_0$ satisfying $\nabla_u (\phi^\oplus(x)+ \gamma(x)\cdot (x_0-x)) \ge 0$ at $x=x_0$, a contradiction. Notice we have shown that if $\supp \pi_{x_0} \neq \{x_0\}$ and a {\sf PDM} exists, then for $c(x,y) = - ||x-y||$, $x_0 \notin \supp \pi_{x_0}$ whenever $ \nabla \phi^\oplus$ exists at $x_0$.

Next, the staying property for the case $c(x,y) =  ||x-y||$ refers to the statement $D_\# (\pi^1 \wedge \pi^2) \leq \pi$ as stated in the theorem. This property was proved in e.g. \cite[Theorem 7.4]{bj} for one dimensional case and then generalized in \cite[Theorem 2.3]{Lim} for general dimension in the {\sf MOT} setup, but note that a solution $\pi$ of {\sf MMOT} is also a solution of {\sf MOT} between its own marginals $\pi^1$ and $\pi^2$. Both \cite{bj}, \cite{Lim} assumed $c(x,y) =|x-y|$, but we note the same proof works for any strictly convex norm cost. From the staying property, the study of the geometry of $\pi$ is now reduced to the study of $\tilde \pi := \pi - D_\# (\pi^1 \wedge \pi^2)$. Note that  $\tilde \pi$ has no mass on the diagonal $\{(x,x) \,|\, x \in \R^d\}$ and $\tilde \pi$ must solve the {\sf VMOT} problem with respect to its own one-dimensional marginals.
\\

{\bf Step 3.} By Step 2 we can assume $\pi_{x_0}$ has no mass at $x_0$. Under this assumption we will show $\supp \pi_{x_0}$ is contained in the set of extreme points of $\conv (\supp \pi_{x_0})$. To this end, first we will show that $\supp \pi_{x_0}$ is contained in the boundary of $\conv (\supp \pi_{x_0})$, where the boundary refers to the topology of $V_0$ and not of $\R^d$. Suppose on the contrary  $\supp \pi_{x_0} \nsubseteq {\rm bd} \big{(}\conv (\supp \pi_{x_0}  )\big{)}$. Then we can find a point $y \in {\rm int} \big{(} \conv (\supp \pi_{x_0}) \big{)} \cap (\supp \pi_{x_0} \setminus \{x_0\})$ and  a subset $\{y_0, y_1,...,y_m\} \subset \supp \pi_{x_0} \setminus \{x_0, y\}$ such that $y$ is a convex combination of these, i.e.,
\begin{align}\label{combi}
\text{There exists } p_i > 0 \text{ for all } i, \sum_{i=0}^m p_i =1 \text{ such that } y= \sum_{i=0}^m p_i y_i.
\end{align}
Then due to the linearity of $y \mapsto  \nabla_u \gamma(x) \cdot (y-x)$,  from \eqref{nablazero}  we deduce
\begin{align}\label{nablaeq}
 \nabla_u c(x_0,y) = \sum_{i=0}^m p_i  \nabla_u c(x_0,y_i) \quad \forall u \in V_0.
  \end{align}
 Using this, we will show that 
 \begin{align}\label{ray}
 \text{the points $\{y, y_0, y_1,...,y_m\}$ lie on a ray emanating from $x_0$.}
 \end{align}
To see this, take $u = y-x_0$, let $c(x,y) = ||x-y||$, and let $\overline{x_0y}$ be the infinite  line  containing $x_0, y$. Thanks to the strict convexity of $|| \cdot ||$, we will firstly  show that $ \nabla_u c(x_0,y) < \nabla_u c(x_0,z)$ for any $z \notin \overline{x_0y}$, where $\nabla_u$ refers to the $u$-directional derivative w.r.t. the $x$ variable. To see this, we consider the following function on $]-1,1[$
 \begin{align*}
\sigma(t) :=  \big{(} ||z-(x_0+tu)|| - ||z-x_0|| \big{)} - \big{(}||y-(x_0+tu) - ||y-x_0|| \big{)}.
  \end{align*}
The claim reads $\sigma ' (0) >0$. Note $||y-x_0||  - ||y-(x_0+tu)|| = t||u||$, so
 \begin{align*}
\sigma(t) =  ||z-(x_0+tu)||  - ||z-x_0|| +  t||u||, \quad -1  < t < 1.
  \end{align*}
Notice $\sigma(0)=0$, $\sigma(t)> 0$ for $t>0$ and $\sigma(t)<0$ for $t<0$ since the norm $|| \cdot ||$ is strictly convex and $z \notin \overline{x_0y}$. 
Also notice $\sigma$ is convex, hence $\sigma ' (0) >0$, showing the claim in this case. Now if $z \in \overline{x_0y}$ with $z$ and $y$ not on the same side of $x_0$, then $\sigma(t) = 2t ||u||$ for small $t$, so that again $\sigma'(0) > 0$ and $ \nabla_u c(x_0,y) < \nabla_u c(x_0,z)$. Finally if $z \in \overline{x_0y}$ with $z$ and $y$ on the same side of $x_0$, then $\sigma'(0) = 0$ and $ \nabla_u c(x_0,y) = \nabla_u c(x_0,z)$. With \eqref{nablaeq} this clearly yields \eqref{ray}. Finally, a similar argument yields \eqref{ray} for the cost $c(x,y) = -||x-y||$.

Now if dim $V_0 \ge 2$, we can choose $\{y_0, y_1,...,y_m\}$  such that they are not aligned meanwhile \eqref{combi} holds. But \eqref{ray} then forces them to be aligned, a contradiction. If dim $V_0 =1$ then since $\pi_{x_0}$ is centered at $x_0$, we can choose $y_0, y_1$ in the opposite direction with respect to $x_0$, i.e. $(y_0 - x_0) \cdot (y_1 -x_0) < 0$. Also $(y_0 - y) \cdot (y_1 - y) < 0$ by \eqref{combi}. Then again \eqref{ray} cannot hold, a contradiction. This yields $\supp \pi_{x_0} \subseteq {\rm bd} \big{(}\conv (\supp \pi_{x_0})\big{)}$. Now if $\supp \pi_{x_0} \nsubseteq {\rm Ext} (\conv (\supp \pi_{x_0}))$, then again we can find $\{y, y_0, y_1,...,y_m\} $ $\subset \supp \pi_{x_0}$ such that \eqref{combi} holds. Then by the above argument we deduce \eqref{ray}, a contradiction to the fact that $\{y, y_0, y_1,...,y_m\} \subseteq {\rm bd} \big{(}\conv (\supp \pi_{x_0})\big{)}$, since a ray emanating from an interior point of a convex set can intersect with its boundary at one point only. Hence we get $\supp \pi_{x_0} \subseteq {\rm Ext} \big{(}\conv (\supp \pi_{x_0})\big{)}$, and as the reverse inclusion is clear, we conclude $\supp \pi_{x_0} = {\rm Ext} \big{(}\conv (\supp \pi_{x_0})\big{)}$. Notice this proves the theorem for $\pi_n = \pi  \big{|}_{ I^n \times \R^d}$, the restriction of $\pi$ on $I^n \times \R^d$ in Step 1. Finally,  letting $n \to\infty$ for the domain $I^n$ completes the proof of Theorem \ref{extreme} for any irreducible pair of marginals $(\mu_i,\nu_i)_i$.
\\

{\bf Step 4.} We will prove the theorem for each pair of marginals $(\mu_i,\nu_i)$ only in convex order and not necessarily irreducible. It is well known that any convex-ordered pair $(\mu_i,\nu_i)$ can be  decomposed as at most countably many irreducible pairs, and the decomposition is uniquely determined by the potential functions $u_{\mu_i}, u_{\nu_i}$. While more details can be found in \cite{bj}, \cite{bnt}, we provide the statement for  reader's convenience.\\

\noindent\cite[Proposition 2.3]{bnt} \,\,For each $i \in [d]$,  let $(I_{i,k})_{1\leq k \leq N}$ be the open components of the open set $\{u_{\mu_i}<u_{\nu_i}\}$ in $\R$, where $N\in \N \cup \{\infty\}$. Let $I_{i,0}=\R\setminus \cup_{k\geq1} I_{i,k}$ and $\mu_{i, k}=\mu  \big{|}_{_{I_{i,k}}}$ for $k\geq 0$, so that $\mu_i=\sum_{k\geq0} \mu_{i,k}$.
  Then there exists a unique decomposition $\nu_i=\sum_{k\geq0} \nu_{i,k}$ such that
  \begin{align*}
  \mu_{i,0} = \nu_{i,0}, \ \mbox{and} \ (\mu_{i,k}, \nu_{i,k}) \mbox{ is irreducible for $k \ge 1$ with } \, \mu_{i,k}(I_{i,k})=\mu_{i,k}(\R).
  \end{align*}
  Moreover, any $\pi_i \in {\rm MT}(\mu_i,\nu_i)$ admits a unique decomposition
  $
    \pi_i=\sum_{k\geq0} \pi_{i,k}
  $
  such that $\pi_{i,k}\in {\rm MT}(\mu_{i,k},\nu_{i,k})$ for all $k\geq0$. \\

Note that $\pi_{i,0}$ must be the identity transport (i.e. $\pi_{i,0}$ is concentrated on the diagonal $\{(x,x) \,|\, x \in I_{i,0} \}$) since it is a martingale and $ \mu_{i,0} = \nu_{i,0}$. There is no randomness in $\pi_{i,0}$. Now for the rest of the proof we will  assume $d=2$ to avoid notational difficulty, but one may observe the same argument works for any $d\ge 3$ just with some notational challenge.

Let $\pi$ be a minimizer in \eqref{VMOT} and let $\mu_1=\sum_{k\geq0} \mu_{1,k}$, $\mu_2=\sum_{k\geq0} \mu_{2,k}$ be the unique decomposition of $\mu_1, \mu_2$ respectively. Then our domain $\R^2$ is decomposed as $\R^2 = \cup_{m \ge 0, n \ge 0} I_{1,m} \times I_{2,n}$ accordingly. Now the strategy is to study the geometry of $\pi$ on each domain $I_{m,n} := I_{1,m} \times I_{2,n}$. 

Let $\pi_{m,n}:= \pi  \big{|}_{ I_{m,n} \times \R^2}$  be the restriction of $\pi$ on $I_{m,n} \times \R^2$. If $m \ge 1$ and  $n \ge 1$, then $(\mu_{1,m}, \nu_{1,m})$ and $(\mu_{2,n}, \nu_{2,n})$ are both irreducible and in this case we already established the theorem. On the other extreme, that is if $m =0$ and  $n =0$, then as mentioned above  there is no randomness in $\pi_{0,0}$, that is $\pi_{0,0} \in {\rm MT}(\pi_{0,0}^1, \pi_{0,0}^1)$ where $\pi_{0,0}^1$ is a coupling of $\mu_{1,0}$ and $\mu_{2,0}$. Hence  $\pi_{0,0}$ represents an identity transport, and thus the theorem obviously holds in this case.

The remaining case is $m=0$ and $n \ge 1$. Write $\pi = \pi_{0,n}$ for simplicity and recall that $\pi = {\rm Law}(X,Y)$ where $X=(X_1,X_2)$, $Y=(Y_1,Y_2)$, and $(X_1,Y_1)$, $(X_2,Y_2)$ are jointly martingales.  Let $\ga_1 = {\rm Law}(X_1,Y_1)$ and $\ga_2 = {\rm Law}(X_2,Y_2)$, so that $\pi$ is a coupling of martingales $\ga_1$ and $\ga_2$. But since $\mu_{1,0}=\nu_{1,0}$ we have $X_1=Y_1$, thus the cost is simplified as follows:
$$c(x,y) = ||(x_1,x_2) - (y_1,y_2)|| = ||(0, x_2 - y_2)|| =: ||x_2 - y_2||_*$$ 
where $|| \cdot ||_*$ is the restriction of the norm $||\cdot||$ on the second coordinate axis in $\R^2$. Hence we have $\E_\pi ||X-Y|| = \E_{\ga_2} ||X_2 - Y_2||_*$ and this implies $\ga_2$ is an optimizer in VMT$(\mu_2, \nu_2)$ with respect to $|| \cdot ||_*$. Now by the irreducibility of $(\mu_2, \nu_2)$ Step 3 already established the theorem for $\ga_2$, hence the theorem also holds for $\pi$ as $\ga_1$ is merely  an identity transport. This completes the proof of the theorem.
\end{proof}

We end this section with some examples that illustrate Theorem \ref{extreme}.

\begin{example}\label{ex1} We give a maximizer for \eqref{VMOT} with $c(x,y)=|x-y|$. As \cite{GKL2} observed, the inequality $\frac{1}{2} (|x-y| - 1)^2 \geq 0$ may be rewritten as
\begin{align}\label{u}
\sum_{i=1}^d \bigg{(}   \frac{1}{2}\big(x_i^2 -\frac1{2d}\big) - \frac{1}{2}\big( y_i^2+\frac1{2d}\big) + x_i (y_i-x_i) \bigg{)} \leq  -|x-y|.
\end{align}
Fix any probability measure $\pi^1 \in \cal P_2(\R^d)$. Choose a family of probabilities $(\pi_x)_x$ on $\R^d$ such that $\delta_x \le_c \pi_x$ and $\pi_x(S_{x,1}) =1$ where $S_{x,1}$ is the unit sphere in $\R^d$ centered at $x$. Now define a martingale measure $\pi$ by $\pi := \pi_x \otimes \pi^1$. Then $\pi$ is a maximizer in ${\rm VMT}(\vec\mu, \vec\nu)$, where $\vec\mu = (\mu_1,...,\mu_d)$ and $\vec\nu = (\nu_1,...,\nu_d)$ are the one-dimensional marginals induced by $\pi$. The optimality of $\pi$ follows from the fact that while the integration of the left hand side of \eqref{u} by any $\gamma \in {\rm VMT}(\vec\mu, \vec\nu)$ is a constant, the inequality \eqref{u} becomes an equality precisely on the set $G:=\{(x,y) :  |x-y|=1\}$, and moreover $\pi(G)=1$.
\end{example}

\begin{example}\label{ex2} We shall construct a minimizer for \eqref{VMOT} with $c(x,y)=|x-y|$. Denote by ${\rm Leb}\big{|}_I$ the Lebesgue measure restricted on a set $I$. 
Let $d=2$, $\mu_1 = \mu_2 = {\rm Leb}\big{|}_{[-1/2, 1/2]}$, $\nu_1= \nu_2 = \frac{1}{2}{\rm Leb}\big{|}_{[-1, 1]}$. Recall that the projection of any element in ${\rm VMT}(\vec\mu, \vec\nu)$ on the ``first martingale space" $\R^2$ belongs to ${\rm MT}(\mu_1, \nu_1)$, that is, if $\E[Y|X]=X$ then $\E[Y_1 | X_1] = X_1$. This implies $P(c) \ge P_1(c)$, where $P_1(c)$ is the primal value of the one-dimensional {\rm MOT} problem
\begin{align}\label{v}
P_1(c) :=\inf_{\pi_1 \in {\rm MT}(\mu_1,\nu_1)}  \int_{\R \times \R} |x_1-y_1|\,d\pi_1(x_1,y_1)
\end{align}
simply because $|x-y| \ge |x_1-y_1|$.
It is well known (see e.g. \cite[Theorem 7.4]{bj}) that the solution to \eqref{v}, say $\ga^1$, is unique and there exist functions $T^- : [-\frac{1}{2}, \frac{1}{2}] \to [-1, -\frac{1}{2}]$, $T^+ : [-\frac{1}{2}, \frac{1}{2}] \to [\frac{1}{2}, 1]$ such that  
\begin{align*}
\ga^1_{x_1} = \frac{1}{2} \bigg{(}\delta_{x_1} +\lambda^-(x_1) \delta_{T^-(x_1)} +\lambda^+(x_1) \delta_{T^+(x_1)}\bigg{)} \ \text{ where} \\
\ga^1 = \ga^1_{x_1} \otimes \mu_1 \ \text{ and } \ \lambda^{\pm} (x_1)=\bigg{|} \frac{T^{\mp} (x_1)-x_1}{T^+(x_1)-T^-(x_1)} \bigg{|}.
\end{align*}
This means that the martingale $\ga^1$ splits the mass at each $x_1 \in ]-1,1[$ in $\mu_1$ onto three points $\{T^-(x_1), x_1, T^+(x_1)\}$. Since the identity transport $x_1 \mapsto x_1$ has no contribution to the cost \eqref{v}, we see that
\begin{align}\label{w}
P_1(c) = \frac{1}{2}\int \big{(} \lambda^-(x_1) |T^-(x_1) - x_1| +  \lambda^+(x_1)   |T^+(x_1) - x_1|  \big{)} d\mu_1(x_1). 
\end{align}
Now let us construct a minimizer $\pi \in {\rm VMT}(\vec\mu, \vec\nu)$. For each $x=(x_1,x_2)\in ]-\frac{1}{2}, \frac{1}{2}[^2$ define a kernel $\pi_x \in \cP(\R^2)$ by
\begin{align*}
\pi_{x} = \frac{1}{2} \bigg{(}&\lambda^-(x_1) \delta_{\big{(}T^-(x_1), x_2\big{)}}  +\lambda^+(x_1) \delta_{\big{(}T^+(x_1), x_2\big{)}} \\
&+\lambda^-(x_2) \delta_{\big{(}x_1, T^-(x_2)\big{)}}    +\lambda^+(x_2) \delta_{\big{(}x_1, T^+(x_2)\big{)}}  \bigg{)}. 
\end{align*}
Note that $\pi_x$ splits mass at $x$ toward four directions -- west, east, south and north. Now let $\pi^1 \in \Pi(\mu_1,\mu_2)$ be any coupling of $\mu_1,\mu_2$, and define $\pi := \pi_x \otimes \pi^1$. Then notice that $\pi \in {\rm VMT}(\vec\mu, \vec\nu)$, and moreover by  \eqref{w}, we have $\int |x-y| d\pi = P_1(c)$. Therefore, $\pi$ is a minimizer.

We note that this example gives a counterexample to the second conjecture given in \cite{GKL2} as $\pi_x$ is supported at four points rather than three. But the existence question of such polytope-type {\rm MOT}s  still remains open. For related results, we refer to \cite{d18-1}.
\end{example}

The next example shows the strict convexity of the norm assumption cannot be omitted in Theorem \ref{extreme}.

\begin{example}
Strict convexity of the norm is necessary for the extremal structure of {\rm VMOT}. For an example, let $d=2$, $c(x,y) = \max(|x_1-y_1|, |x_2-y_2|)$, $\mu_1 = \mu_2 = {\rm Leb}\big{|}_{[-1/2, 1/2]}$, $\nu_1 = \frac{1}{2}(\delta_{-10} + \delta_{10})$, $ \nu_2 = \frac{1}{2}{\rm Leb}\big{|}_{[-1, 1]}$. Observe every $\gamma \in {\rm VMT}(\vec\mu, \vec\nu)$ leads to the same cost $\int |x_1-y_1|\,d\ga^1$ where $\gamma^1 \in {\rm MT}(\mu_1,\nu_1)$ (note ${\rm MT}(\mu_1,\nu_1)$ is a singleton). Now clearly there are elements in ${\rm VMT}(\vec\mu, \vec\nu)$ which do not satisfy Theorem \ref{extreme}, for example one may take $\ga = \ga^1 \otimes \ga^2$ where $\ga^2 \in {\rm MT}(\mu_2,\nu_2)$ is independent of $\ga^1$ and its kernel $\{\ga^2_x\}_x$ are dispersed onto  $[-1,1]$.
\end{example}

\section{Existence of dual optimizers for VMOT -- Theorem \ref{main2}}\label{m2}

One of the most important results for the proof of Theorem \ref{main2} is Proposition \ref{conv}, a compactness property of a class of convex functions. 

Recall two probabilities in convex order 
$\mu\leq_{c}\nu$ is irreducible if $I =\{x \, | \, u_{\mu}(x) < u_{\nu}(x)\}$ is connected and $\mu(I)=\mu(\R)$; see \cite[Section 2]{bnt}. In this case, $(I,J)$ is called the \emph{domain} of $(\mu,\nu)$ where $J$ is the smallest interval satisfying $\nu(J) = \nu(\R)$, that is, $J$ is the union of~$I$ and any endpoints of $I$ that are atoms of $\nu$, and moreover, $I={\rm int}(J) = {\rm int}({\rm conv}(\supp (\nu)))$. Note that  $\mu\leq_{c}\nu$ is irreducible if and only if for every $\pi  =\pi_x \otimes \mu \in {\rm MT}(\mu,\nu)$ and $z \in I:={\rm int}(\conv(\supp \nu))$, we have
\begin{align}\label{irreducibility}
\text{$\mu(I^\pi_z) >0$ \quad where \quad $I^\pi_z:= \{ x \in I \,|\, z \in {\rm int}({\rm conv}(\supp (\pi_x))) \}$.}
\end{align}
In other words, $(\mu,\nu)$ is irreducible if and only if for any $z \in {\rm int(conv}(\supp \nu))$, every $\pi \in {\rm MT}(\mu,\nu)$ must ``cut across" the point $z$.
This follows directly from the definition of potential functions and Jensen's inequality.

For  irreducible pairs $(\mu_i,\nu_i)_{i \in [d]}$ with domains $(I_i, J_i)$, let us define $I = I_1 \times ... \times I_d$, $J = J_1 \times ... \times  J_d$ and  $\mu^\otimes = \mu_1 \otimes ... \otimes \mu_d$,  $\nu^\otimes = \nu_1 \otimes ... \otimes \nu_d$. 

\begin{proposition}\label{conv} Let $(\mu_i,\nu_i)_{i \in [d]}$ be  irreducible pairs of probabilities with domains $(I_i,J_i)_i$. Let $a \in I$, $C \in \R$. Consider the following class of functions $\Lambda=\Lambda(a,C,\vec\mu,\vec\nu) $ where every $\chi \in \Lambda$ satisfies the following:
\begin{enumerate}
\item $\chi$ is a real-valued convex function on $J$,
\item $\chi \ge 0$ and $\chi(a)=0$,
\item $\int \chi \,d(\nu^\otimes - \mu^\otimes) \le C$.
\end{enumerate}
Then $\Lambda$ is locally bounded in the following sense: for each compact subset $K$ of $J$, there exists $M=M(K)$ such that $\chi \le M$ on $K$ for every $\chi \in \Lambda$. Moreover, for any sequence $\{\chi_n\}_n$ in $\Lambda$ there exists a subsequence $\{ \chi_{n_j} \}_j$ of $\{ \chi_n\}_n$ and a real-valued convex function $\chi$ on $J$ such that $\lim_{j \to \infty} \chi_{n_j} (x) = \chi(x)$ for every $x \in J$, and  the convergence is uniform on every compact subset of $I$.
\end{proposition}
We note that the proposition for the case $d=1$ was proved in \cite{bnt}. In this section we prove Theorem \ref{main2} while assuming Proposition \ref{conv}, whose proof will then be given in an appendix. During the proof the bounding constant $C$ may vary, but observe it does not depend on $n$.
\begin{proof} 
{\bf Step 1.}  The assumption $|c(x,y)| \le \sum_{i=1}^d \big{(}v_i(x_i) + w_i(y_i)\big{)}$ for continuous $v_i \in L^1(\mu_i)$, $w_i \in L^1(\nu_i)$ ensures $P(c)=D(c)$ as in \eqref{duality} (for a proof, see e.g. \cite{Z}). Moreover, it is obvious that a dual optimizer exists for $c(x,y)$ iff so does for $\tilde c (x,y):= c(x,y) - \sum_{i=1}^d \big{(}v_i(x_i) + w_i(y_i)\big{)}$. Hence by replacing $c$ with $\tilde c$, from now on we will assume that $c \le 0$.

As $P(c)=D(c) \in \R$, we can find an ``approximating dual maximizer" $(f_{i,n}, g_{i,n}, h_{i,n})_{n \in \N}$ which consist of  real-valued continuous functions $f_{i,n} \in L^1(\mu_i)$, $g_{i,n} \in L^1(\nu_i)$, and continuous bounded $h_{i,n}$ for every $i \in [d]$ and $n \in \N$,  such that the following weak duality holds:
\begin{align}
\label{dual} &\sum_{i=1}^d \big{(} f_{i,n} (x_i) - g_{i,n} (y_i) + h_{i,n}(x) (y_i - x_i) \big{)} 
 \le c(x,y),\\
\label{maximizing} &\sum_{i=1}^d  \big{( } \int f_{i,n}d\mu_i - \int g_{i,n}d\nu_i \big{)}  \nearrow P(c)
\quad \text{as} \quad  n \to \infty.
\end{align}
Denote $f_n^\oplus(x) =\sum_{i=1}^d  f_{i,n} (x_i)$, $g_n^\oplus(y) = \sum_{i=1}^d g_{i,n} (y_i)$, $h_n(x) = \big{(} h_{1,n}(x),..., h_{d,n}(x)\big{)}$. By taking supremum over $x \in \R^d$ in \eqref{dual}, we get
\begin{align}
\label{chidefinition}\chi_n(y):= \sup_{x \in \R^d} \big{(} f_n^\oplus (x) + h_{n}(x) \cdot (y - x) \big{)} \le g_n^\oplus(y).
\end{align}
Notice $\chi_n$ is a convex function on $\R^d$, and by definition of $\chi_n$, we see 
\begin{align}\label{ineq}
f_n^\oplus  \le  \chi_n \le g_n^\oplus \quad \forall n \in \N.
\end{align}
In this step we shall obtain local uniform boundedness of $\{\chi_n\}_n$ \eqref{boundconvex}. 
Let $a \in I$ be fixed.  By subtracting an appropriate linear function $L_n(y) = \nabla L_n(x) \cdot (y - x) + L_n(x)$ from \eqref{dual}, that is by replacing $f_n^\oplus(x)$ with $f_n^\oplus(x)- L_n(x)$, $g_n^\oplus(y)$ with $g_n^\oplus(y)- L_n(y)$ and $h_n(x)$ with $h_n(x) - \nabla L_n(x)$, we can assume 
\begin{align}\label{convex0}
\chi_{n} \ge 0 \ \text{ and } \ \chi_{n}(a) =0 \quad \forall n.
\end{align}
Next, observe \eqref{ineq} gives
\[
\sum_i \int g_{i,n} d\nu_i \ge \int \chi_n d\nu^\oplus \ge \int \chi_n d\mu^\oplus \ge \sum_i\int f_{i,n} d\mu_i
\]
where the second inequality is due to the convexity of $\chi_n$ and $\mu_i \le_c \nu_i$. Then by \eqref{maximizing}, there exists a constant $C >0$ such that
\begin{align}\label{intbypart}
\int \chi_{n} \,d(\nu^\otimes - \mu^\otimes) \le C \quad \forall n.
\end{align}
Let $\{\epsilon_k\}_k$ be a positive decreasing sequence tending to zero as $k \to \infty$, and write $I_i = ]a_i, b_i[$ where $-\infty \le a_i < b_i \le +\infty$. Then we define the compact interval $J_{i,k}:=[c_{i,k}, d_{i,k}]$ for $i \in [d]$ and $k \in \N$ as
\begin{align}
&\text{If } a_i > -\infty \text{ then: } \ \nu_i(a_i)=0 \Rightarrow c_{i,k }:= a_i + \epsilon_k, \  \nu_i(a_i) > 0 \Rightarrow c_{i,k} := a_i, 
\nonumber \\
&\text{If } b_i < +\infty \text{ then: } \ \nu_i(b_i)=0 \Rightarrow d_{i,k }:= b_i -  \epsilon_k,\ \nu_i(b_i)>0 \Rightarrow d_{i,k} := b_i,  
\nonumber \\
&\text{If }  a_i = -\infty \Rightarrow c_{i,k} := -1/ \epsilon_k.  \ \text{If }  b_i = +\infty \Rightarrow d_{i,k} := +1/\epsilon_k. \nn
\end{align}
Thus for example, if $\nu(a_i) = 0$ and $\nu(b_i) > 0$, then $J_{i,k} = [a_i + \epsilon_k, b_i]$. Let $\epsilon_1$ be small so that $\mu_i (J_{i,1})>0$, $\nu_i (J_{i,1}) > 0$ for every $i \in [d]$. Notice $J_{i,k} \nearrow  J_i$. Let $J_k := J_{1,k} \times J_{2,k} \times ... \times J_{d,k}$. 
Then by Proposition \ref{conv},  there exists  $\{M_k\}_k$ such that 
\begin{align}\label{boundconvex}
0 \le \sup_n \chi_{n} \le M_k \ \ \text{on} \ \ J_k.
\end{align}

{\bf Step 2.} Let  $(f_{i,n}, g_{i,n}, h_{i,n})_{n \in \N}$ be an approximating dual maximizer. We want to show pointwise convergence of $f_{i,n}$ and $g_{i,n}$, i.e. $ f_{i,n}(x_i)\to f_i(x_i) \in \R$ for $\mu_i$-a.e. $x_i$ and  $ g_{i,n}(y_i) \to g_i(y_i) \in \R$ for $\nu_i$-a.e. $y_i$ as $n \to \infty$. However, in establishing the convergence we see that there is an immediate obstacle when $d \ge 2$, that is, in the duality formulation \eqref{dual} one can always replace $(f_{i,n})_{i \in [d]}$ by $( f_{i,n}+ C_{i,n})_{i \in [d]} $ for any constants $(C_{i,n})_{i \in [d]} $ satisfying 
  $\sum_i C_{i,n} = 0$ and similarly replace $(g_{i,n})_{i \in [d]}$ by $( g_{i,n}+ D_{i,n})_{i \in [d]} $. This means the convergence cannot be shown for any approximating dual maximizer. In this step we will show there exists an approximating dual maximizer which satisfies the convergence. 

Take an approximating dual maximizer $(f_{i,n}, g_{i,n}, h_{i,n})_{n \in \N}$. Observe 
\begin{align*}
C &\ge \int g_n^\oplus d\nu^\otimes - \int f_n^\oplus d\mu^\otimes   \\
&\ge \int \chi_n d\nu^\otimes- \int f_n^\oplus d\mu^\otimes
\\
&\ge \int \chi_n d\mu^\otimes- \int f_n^\oplus d\mu^\otimes  \\
&=\, \, \parallel \chi_n - f_n^\oplus \parallel_{L^1(\mu^\otimes)} \quad  \forall n
\end{align*}
where the third inequality is by convexity of $\chi_n$ and $\mu^\otimes \le_c \nu^\otimes$. For each $k \in \N$, let $\mu_{i,k}$ be the restriction of $\mu_i$ on $J_{i,k}$ then normalized to be a probability measure. Let $\mu^\otimes_k = \otimes_i \mu_{i,k}$, so that $\mu^\otimes_k(J_k)=1$. Define
\begin{align*}
v_{i,k,n} = \int f_{i,n} \,d\mu_{i,k}, \quad i \in [d], k\in \N ,n \in \N.
\end{align*}
Let us fix $k \in\N$. We claim that 
\begin{align*}
\sup_n \parallel f_{i,n} - v_{i,k,n} \parallel_{L^1(\mu_{i,k})} \text{ is bounded for each $i\in[d]$}.
\end{align*}
To see this, recall $ \sup_n \parallel \chi_n - f_n^\oplus \parallel_{L^1(\mu^\otimes_k)}$ is bounded  and  so by \eqref{boundconvex},
\begin{align*}
\parallel M_k - f_n^\oplus \parallel_{L^1(\mu^\otimes_k)} =\parallel M_k - \chi_n \parallel_{L^1(\mu^\otimes_k)} + \parallel \chi_n - f_n^\oplus \parallel_{L^1(\mu^\otimes_k)} \, \le M_k +C =: C.
\end{align*}
The constant $C$ may depend on $k$ but not on $n$. From this and Jensen,
\begin{align}\label{boundv}
|v_{1,k,n} + v_{2,k,n}+...+v_{d,k,n}| \le \, \parallel f_n^\oplus \parallel_{L^1(\mu^\otimes_k)}  \, \le C \, \text{ for all } n.
\end{align}
Next, note that as $f_n^\oplus \le M_k$ on $J_k$, by taking supremum we see that
\begin{align*}
\sum_{i=1}^d \sup_{x_i \in J_{i,k}} f_{i,n} (x_i) \le M_k \q \forall n,
\end{align*}
and note that obviously $v_{i,k,n} \le \sup_{x_i \in J_{i,k}} f_{i,n} (x_i)$, so in particular
\begin{align*}
\sup_{x_1 \in J_{1,k}} f_{1,n} (x_1)+ \sum_{i=2}^d v_{i,k,n}  \le M_k.
\end{align*}
Define  $\hat v_{1,k,n} = -\sum_{i=2}^d v_{i,k,n}$ and observe
\begin{align*}
C \ge \,\,\parallel M_k - f_n^\oplus \parallel_{L^1(\mu^\otimes_k)} = M_k - \int ( f_{1,n} + \sum_{i=2}^d v_{i,k,n}) d\mu_{1,k}.
\end{align*}
This implies that $ \sup_n \parallel  f_{1,n} - \hat v_{1,k,n} \parallel_{L^1(\mu_{1,k})}$ is bounded, and then by \eqref{boundv}, $ \sup_n \parallel  f_{1,n} - v_{1,k,n} \parallel_{L^1(\mu_{1,k})}$ is bounded. The claim is proved.

We are now ready to apply the Koml{\'o}s lemma, which states that every $L^1$-bounded sequence of real functions contains a subsequence such that the arithmetic means of all its subsequences converge pointwise almost everywhere. Define $\tilde v_{1,k,n}=  \frac{1}{n} \sum_{m=1}^n \hat v_{1,k,m}$, and $\tilde v_{i,k,n}=  \frac{1}{n} \sum_{m=1}^n v_{i,k,m}$ for $2 \le i \le d$.
Observe that by repeated use of Koml{\'o}s lemma, for every $i\in[d]$ and $k \in \N$, we can find a subsequence $\{f_{i,k,n}\}_n$ of $\{f_{i,n}\}_n$ such that 
\begin{enumerate}
\item $\{f_{i,k+1,n}\}_n$ is a further susequence of $\{f_{i,k,n}\}_n$, and
\item $\tilde f_{i,k,n}(x_i) - \tilde v_{i,k,n}$ converges for $\mu_{i,k}$ - a.e. $x_i$ as $n \to \infty$,
\end{enumerate}
where $\tilde f_{i,k,n}(x) = \frac{1}{n} \sum_{m=1}^n f_{i,k,m} (x)$. Select the diagonal sequence  $F_{i,n} := f_{i,n,n}$ and again define $w_{i,k,n} = \int F_{i,n} d\mu_{i,k}$ for $2\le i \le d$, $\hat w_{1,k,n} = -\sum_{i=2}^d w_{i,k,n}$, and $\tilde F_{i,n}(x_i) = \frac{1}{n} \sum_{m=1}^n F_{i,m} (x_i)$, $\tilde w_{1,k,n}=  \frac{1}{n} \sum_{m=1}^n \hat w_{1,k,m}$, and $\tilde w_{i,k,n}=  \frac{1}{n} \sum_{m=1}^n w_{i,k,m}$ for $2 \le i \le d$. We finally claim that 
 \begin{align}\label{convergeF}
\tilde F_{i,n}(x_i)  - \tilde w_{i,1,n} \ \text{converges for} \ \mu_i - a.e.\,x_i \ \text{ for every } i \in [d].
\end{align}
The point is that the dependence on $k$ is now removed. To see this, as $\{F_{i,n}\}_n$ is a subsequence of $\{f_{i,k,n}\}_n$ for every $k\in \N$,  by Koml{\'o}s lemma
 \begin{align}\label{convergeFk}
\tilde F_{i,n}(x_i)  - \tilde w_{i,k,n} \ \text{converges for} \ \mu_{i,k} - a.e.\,x_i \ \text{ for every } i \in [d].
\end{align}
In particular, both $\{\tilde F_{i,n}(x_i)  - \tilde w_{i,1,n}\}_n$ and $\{\tilde F_{i,n}(x_i)  - \tilde w_{i,k,n}\}_n$ converge for $ \mu_{i,1}$ - a.e. $x_i$ as $n\to \infty$, hence their difference $ \{\tilde w_{i,1,n} -  \tilde w_{i,k,n}\}_n$ must converge for any fixed $k$. With \eqref{convergeFk} this implies \eqref{convergeF}.  

Now having an approximating dual maximizer $(F_{i,n}, G_{i,n}, H_{i,n})_n$ where $G,H$ followed the same subsequence as $F$ did (and noting that \eqref{dual} and \eqref{maximizing} are conserved while taking convex combinations), we can repeat the same procedure for $G_{i,n}$ by starting from the uniform bound $\parallel G_n^\oplus - \chi_n \parallel_{L^1(\nu^\otimes)} \, \le C$ so that we get a further subsequence for which similar statement as in \eqref{convergeF} holds for $G_{i,n}$ as well. Then a final  application of Koml{\'o}s lemma (recall $\sum_i \tilde w_{i,1,n} = 0$) yields an approximating dual maximizer $(\tilde F_{i,n}, \tilde G_{i,n}, \tilde H_{i,n})_n$ which satisfies the claimed convergence property.\\

{\bf Step 3.} We've shown there exists a sequence of functions $(f_{i,n}, g_{i,n}, h_{n})_n $ satisfying \eqref{dual}, \eqref{maximizing} for each $n$, and admitting the limit functions $(f_i, g_i)_{i \in [d]}$ such that as $n \to \infty$, $f_{i,n}(x_i) \to f_i(x_i)$ for $\mu_i$ - a.e. $x_i$ and $g_{i,n}(y_i) \to g_i(y_i)$ for $\nu_i$ - a.e. $y_i$. For convenience, define $A_i$ and $B_i$ to be the set of convergence points, that is $\mu_i(A_i) =1$, $\nu_i(B_i)=1$, and 
$$\lim_{n \to \infty} f_{i,n}(x_i) = f_i(x_i) \in \R, \  \lim_{n \to \infty} g_{i,n}(y_i) = g_i(y_i) \in \R \quad \forall x_i \in A_i,\, y_i \in B_i.$$
Let $A = A_1 \times...\times A_d$, $B = B_1 \times...\times B_d$. Note that $A \subset I$, $B \subset J$, and moreover the interior of the convex hull of $B$ is $I$, i.e. ${\rm int}({\rm conv}(B))=I$. 

In this step we will show the convergence of $\{\chi_n\}_n$ defined in \eqref{chidefinition}. Fix $a \in A$ so that $\lim_{n \to \infty} f^\oplus_n(a) = f^\oplus (a)=:\sum_i f_i(a)$. As ${\rm int}({\rm conv}(B))=I$, we can find points in $B$, say $\{y_1,...,y_m\}$, such that $\lim_{n \to \infty} g^\oplus_n(y_j) = g^\oplus (y_j)$ for every $1 \le j \le m$ and $a \in {\rm int}({\rm conv}(\{y_1,...,y_m\}))$. Then in view of \eqref{ineq} we see that both $\{\chi_n(a)\}_n$ and $\{\nabla \chi_n (a)\}_n$ are uniformly bounded in $n$, where $\nabla \chi_n (a) \in \partial \chi_n (a)$ is  a subgradient  of the convex function $\chi_n$ at $a$. Hence by taking a subsequence, we can assume that $\{\chi_n(a)\}_n$ and $\{\nabla \chi_n (a)\}_n$ both converge. Now define the affine function $L_n(y) = \chi_n(a) + \nabla \chi_n (a) \cdot (y-a)$ and replace $f_n^\oplus(x)$ with $f_n^\oplus(x)- L_n(x)$, $g_n^\oplus(y)$ with $g_n^\oplus(y)- L_n(y)$ and $h_n(x)$ with $h_n(x) - \nabla \chi_n(a)$. Then we get $\chi_n(a) = \nabla \chi_n(a) = 0$ while the convergence property of $f_{i,n}, g_{i,n}$ is retained. And \eqref{intbypart} and Proposition \ref{conv} yield $\lim_{n \to \infty} \chi_n = \chi$ on $J$.
\\

{\bf Step 4.}  We will show there exists a function $h : I \to \R^d$ such that 
\begin{align}\label{duallimit}
\sum_{i \in [d]}\big{(}f_i(x_i)-g_i(y_i)\big{)} + h(x) \cdot (y-x) \le c(x,y) \q \forall x,y \in \R^d.
\end{align}
We may define $f_i = -\infty$ on $\R \setminus A_i$ and $g_i= +\infty$ on $\R \setminus B_i$ so that they are defined everywhere on $\R$. For a function $\phi$ defined on a subset of $\R^d$ and bounded below by an affine function, let ${\rm conv}[\phi]:\R^d \to \R \cup \{\infty\}$ denote the convex envelope of $\phi$. Define 
\begin{align}\label{H}
H(x,y) := {\rm conv}[c(x,\cdot) + g^\oplus (\cdot)](y).
\end{align}
Then any measurable choice $h$ satisfying $h(x) \in \partial H(x, \cdot)(x)$ will satisfy \eqref{duallimit}. To see this, we may argue similarly as \cite{bnt}. Recall
\begin{align*}
f_n^\oplus(x) +h_n(x) \cdot (y-x) \le c(x,y) + g_n^\oplus(y).\end{align*}
If we define $H_n(x,y) := {\rm conv}[c(x,\cdot) + g_n^\oplus(\cdot)](y)$, then we have 
\begin{align*}
f_n^\oplus(x) +h_n(x) \cdot (y-x) \le H_n(x,y)  \le c(x,y) + g_n^\oplus(y).
\end{align*}
In particular by taking $y=x$, we get $
f_n^\oplus(x) \le H_n(x,x)$. 
Next, observe that since the $\limsup$ of convex functions is convex, we have
\begin{align*}
\limsup_{n \to \infty} H_n(x,y) &\le {\rm conv}[\limsup_{n \to \infty}\big{(} c(x,\cdot) + g_n^\oplus(\cdot)\big{)}](y)  \\
&= {\rm conv}[c(x,\cdot) + g^\oplus(\cdot)](y)= H(x,y).
\end{align*}
Then by the convergence $f_n^\oplus \to f$ and the  definition of $H(x,y)$, we get
\[
f(x) \le H(x,x), \ \text{and} \ H(x,y) \le c(x,y) +g^\oplus(y).
\]
Since ${\rm int}({\rm conv}(B))=I$ and $g^\oplus$ is real-valued on $B$, for each $x\in I$ the convex function $y \mapsto H(x,y)$ is real-valued thus continuous in $I$. The subdifferential $\partial H(x, \cdot) (y)$ is then nonempty, convex and compact for every $y \in I$. Hence we can  choose a measurable function $h : I \to \R^d$ satisfying $h(x) \in \partial H(x, \cdot)(x)$. Any such choice  yields \eqref{duallimit} as follows:
 \[
 f^\oplus(x) + h(x) \cdot (y-x) \le H(x,x) + h(x) \cdot (y-x) \le H(x,y) \le c(x,y) + g^\oplus(y).
  \]
\smallskip

{\bf Step 5.} We will show that for any function $h : I \to \R^d$ satisfying 
\begin{align}\label{pointwisedualineq}
f^\oplus(x)-g^\oplus(y) + h(x) \cdot (y-x) \le c(x,y)
\end{align}
(we know that such a function exists by Step 4), and for any minimizer $\pi^*$ for the problem \eqref{VMOT}, in fact we have
\begin{align}\label{pointwisedualeq}
f^\oplus(x) - g^\oplus(y)+ h(x) \cdot (y-x) = c(x,y), \quad \pi^* - a.s..
\end{align}
In other words, every minimizer $\pi^*$ is concentrated on the contact set
$$\Gamma := \{(x,y) \,|\, f^\oplus(x) - g^\oplus(y)+ h(x) \cdot (y-x) = c(x,y) \}$$
whenever $h$ is chosen to satisfy \eqref{pointwisedualineq}. This will complete the proof.

To begin, recall $f_{i,n} \to f_i$ $\mu_i$-a.s., $g_{i,n} \to g_i$ $\nu_i$-a.s., and $\chi_n \to \chi$ on $J$ where $\chi_n$ is defined as \eqref{chidefinition} with $\chi$ being its limit. For any $\pi \in {\rm VMT}(\vec\mu,\vec\nu)$ (not necessarily an optimizer), notice that $c(x,y) \in L^1(\pi)$ by the assumption of Theorem \ref{main2}. Now we claim:
\begin{align*}
\limsup_{n \to \infty}& \int \big{(}f_n^\oplus(x)-g_n^\oplus(y) + h_n(x) \cdot (y-x) \big{)}  d\pi  \\
\le & \int \big{(}f^\oplus(x)-g^\oplus(y) + h(x) \cdot (y-x) \big{)}  d\pi.
\end{align*}
First, let us observe how the claim implies \eqref{pointwisedualeq}. Let $\pi^*$ be any minimizer for \eqref{VMOT} (note that this exists by the assumption on the cost). Then $c(x,y) \in L^1(\pi^*)$ and $P(c) = \int c(x,y) \,d\pi^*$. Now we deduce
\begin{align*}
P(c)&=\lim_{n \to \infty}  \int \big{(}f_n^\oplus(x)-g_n^\oplus(y) + h_n(x) \cdot (y-x) \big{)}  d\pi^* \\ 
&\le \int \big{(}f^\oplus(x)-g^\oplus(y) + h(x) \cdot (y-x) \big{)}  d\pi^* \\ 
&\le  \int c(x,y) \,d\pi^* \\  
& = P(c)
\end{align*}
hence equality holds throughout, and this implies \eqref{pointwisedualeq}.

If we have pointwise convergence $h_n(x) \to h(x)$ then 
the claim would have been a simple consequence of Fatou's lemma, but we do not know such a convergence a priori. But \cite{bnt} suggested a clever idea to handle this case and let us carry out a similar scheme in the current $\R^d$ setting. 

Fix any $\pi \in {\rm VMT}(\vec\mu,\vec\nu)$ and let $\pi^1 := p^x_\# \pi$, $\pi^2 := p^y_\# \pi$. Then as in Step 2 (but $\pi^1 \le_c \pi^2$ instead of $\mu^\otimes \le_c \nu^\otimes$), by \eqref{maximizing} and \eqref{ineq}  we get
\begin{align*}
\sup_n || \chi_n - f_n^\oplus ||_{L^1(\pi^1)} < \infty, \quad \sup_n ||g_n^\oplus -  \chi_n  ||_{L^1(\pi^2)} < \infty.
\end{align*}
From this, as $f_n^\oplus \to f^\oplus$, $g_n^\oplus \to g^\oplus$, $\chi_n \to \chi$, by Fatou's lemma we get
\begin{align*}
&\chi - f^\oplus \in L^1(\pi^1), \quad g^\oplus -  \chi  \in L^1(\pi^2),\\
&\limsup_{n \to \infty} \int (f_n^\oplus - \chi_n )\,d\pi^1 \le  \int (f^\oplus - \chi )\,d\pi^1,\\
 &\limsup_{n \to \infty} \int ( \chi_n -g_n^\oplus)\,d\pi^2 \le  \int ( \chi - g^\oplus )\,d\pi^2.
\end{align*}
This allows us to proceed
\begin{align*}
\limsup_{n \to \infty}  &\int \big{(}f_n^\oplus(x)-g_n^\oplus(y) + h_n(x) \cdot (y-x) \big{)}  d\pi   \\
= \limsup_{n \to \infty}  &\int \big{(}f_n^\oplus(x)- \chi_n(x) + \chi_n(y) -g_n^\oplus(y) \\
&+\chi_n(x) - \chi_n(y) + h_n(x) \cdot (y-x) \big{)} d\pi  \\
\le &\int (f^\oplus - \chi )\,d\pi^1 + \int ( \chi - g^\oplus )\,d\pi^2 \\ &+ \limsup_{n \to \infty} \int \big{(}\chi_n(x) - \chi_n(y) + h_n(x) \cdot (y-x) \big{)} \, d\pi.
\end{align*}
To handle the last term, let $\pi = \pi_x \otimes \pi^1$, and let $\xi_n : I \to \R^d$ be a sequence of functions satisfying $\xi_n(x) \in \partial \chi_n(x)$. Then we compute
\begin{align*}
&\int \big{(}\chi_n(x) - \chi_n(y) + h_n(x) \cdot (y-x) \big{)} \, d\pi   \\
=& \iint  \big{(} \chi_n(x) - \chi_n(y) + h_n(x) \cdot (y-x) \big{)} \, d\pi_x (y) \,d\pi^1(x)   \\
=&\iint  \big{(} \chi_n(x) - \chi_n(y) +\xi_n(x) \cdot (y-x) \big{)} \, d\pi_x (y) \,d\pi^1(x),
\end{align*}
since $\int h_n(x) \cdot (y-x) \, d\pi_x (y)$ = $\int \xi_n(x) \cdot (y-x) \, d\pi_x (y)=0$. Notice that the last integrand is nonpositive. Thus by Fatou's lemma,
\begin{align*}
& \limsup_{n \to \infty} \int \big{(}\chi_n(x) - \chi_n(y) + h_n(x) \cdot (y-x) \big{)} \, d\pi
   \\
&\le \int \limsup_{n \to \infty}\bigg{(} \int  \big{(} \chi_n(x) - \chi_n(y) + \xi_n(x) \cdot (y-x) \big{)} \, d\pi_x (y)\bigg{)} \,d\pi^1(x)   \\
&\le \int \bigg{(} \int \big{(} \chi(x) - \chi(y) + \xi(x) \cdot (y-x) \big{)} \, d\pi_x (y) \bigg{)} \,d\pi^1(x)   \\
\end{align*}
for some $\xi(x)  \in \partial \chi(x)$ which is a limit point of the bounded sequence $\{ \xi_n(x)\}_n$. Finally, in the last line, the inner integral equals
\[\int \big{(} \chi(x) + h(x) \cdot (y-x) - \chi(y) \big{)} \, d\pi_x (y).\]
This proves the claim, hence the theorem.
\end{proof}

\begin{remark}\label{remark3.2}
In Step 5 we deduced $\chi - f^\oplus \in L^1(\pi^1)$, \,$g^\oplus -  \chi  \in L^1(\pi^2)$ for every $\pi \in {\rm VMT}(\vec\mu,\vec\nu)$. Is $f_i \in L^1(\mu_i)$ or $g_i \in L^1(\nu_i)$? Even in one-dimension and in a fairly mild situation such as $c(x,y)$ is $1$-Lipschitz and $\mu, \nu$ are irreducible and compactly supported, it can happen that neither $f \in L^1(\mu)$ nor $g \in L^1(\nu)$, as shown by examples in \cite{bnt}. We may extend the notion of generalized integral of the pair $(f,g)$ introduced in \cite{bnt} to the {\sf VMOT} setup via (compare with \cite[Definition 4.7]{bnt})
 \begin{align}
 \pi^1(f^\oplus) -  \pi^2(g^\oplus) := \pi^1(f^\oplus - \chi) - \pi^2(g^\oplus - \chi) + (\pi^1 - \pi^2) (\chi).
 \end{align}
 For this definition to be meaningful, it is desired that the right hand side would not depend on the choice of $\pi$, since this will then allow us to write the left hand side as e.g. $\vec \mu(f^\oplus) - \vec\nu(g^\oplus)$. Although this looks fairly plausible, we shall not pursue this point further in this paper.
\end{remark}

\section{Further results on the structure of VMOT}\label{m3}
In this section, we present two more structural results for {\sf VMOT}. Recall the proof of Theorem \ref{extreme} was based on the differential identity
\[
\nabla_x \bigg{(}\sum_{i \in [d]}\big{(} f_i(x_i)-g_i(y_i) + h_i(x)  (y_i-x_i)\big{)}  - c(x,y)\bigg{)} = 0, \ \ \pi^* - a.e.\, (x,y)
\]
where $\pi^*$ is a {\sf VMOT} and $(f_i,-g_i,h_i)_i$ is a {\sf PDM}. On the other hand, the next result shall be obtained via the differential identity in $y$ variable:
\[
\nabla_y \bigg{(}\sum_{i \in [d]}\big{(} f_i(x_i)-g_i(y_i) + h_i(x) (y_i-x_i)\big{)}  - c(x,y)\bigg{)} = 0, \ \ \pi^* - a.e.\, (x,y).
\]

\begin{theorem}\label{mongesolution}
Assume the same as in Theorem \ref{main2}. Let $S$ be a proper subset of $[d]$ and assume that $\nu_i$ is absolutely continuous with respect to  Lebesgue measure for every $i \in S$. Suppose $c(x,y)$ is semiconcave in $y$ in the following sense: there exists $u_i : J_i \to \R$ for all $i \in [d]$ such that
\begin{align*}
y \mapsto c(x,y) + \sum_i u_i(y_i) \text{ is concave on $J$ for every $x\in I$.} 
\end{align*}
Assume the following twist condition: for each $x \in I$ and $y_i \in J_i$ where $i \in S$, the mapping
 \begin{align}\label{twist}
(y_j)_{j \in [d] \setminus S} \mapsto \bigg{(} \frac{\partial c(x, y)}{\partial y_i}  \bigg{)}_{i \in S} \text{ is one-to-one on }  \underset{j \in [d] \setminus S}{\otimes} J_j 
\end{align}
whenever the derivatives exist. Then there exists a family of functions $F_x : \R^{|S|} \to \R^{d-|S|}$ for each $x \in \R^d$ such that for any minimizer $\pi = \pi_x \otimes \pi^1$ of \eqref{VMOT}, $\pi_x$ is concentrated on the graph of $F_x$ for $\pi^1$-a.e. $x$.
\end{theorem}
 
We note that the definition of semiconcavity here is so broad that there is no regularity imposed on $u_i$ other than that it be real-valued and measurable.

\begin{proof}
The same assumption as in Theorem \ref{main2} ensures there is a {\sf PDM}. By the semiconcavity on the cost, we can assume that $y \mapsto c(x,y)$ is concave for every $x \in I$ for which $(f_i, -g_i, h_i)_i$ is a {\sf PDM}, that is
\begin{align*}
f^\oplus(x) + h(x) \cdot (y-x) - c(x,y) \le g^\oplus(y) \q \forall x,y \in \R^d,
\end{align*}
and for any minimizer $\pi\in {\rm VMT}(\vec\mu,\vec\nu)$ of the primal problem \eqref{VMOT}, 
\begin{align*}
f^\oplus(x) + h(x) \cdot (y-x) - c(x,y) = g^\oplus(y)  \quad \pi - a.e.\, (x,y).
\end{align*}
Consider the ``inverse martingale Legendre transform" (see \cite{GKL2})
\begin{align}\label{1beta}
\beta(y) := \sup_{x\in I } \{ f^\oplus(x) + h(x) \cdot (y-x) - c(x,  y)\}.
\end{align}
Notice $\beta$ is convex and $\beta \le \sum_i g_i$. By replacing $g_i$'s with its Legendre transform $\psi_i$'s with respect to the ``cost" $\beta$, that is, defining $\psi_i$ successively for $i=1,2,...,d$ by
\[
\psi_i(y_i) := \sup_{y_j,\, j \neq i} \big{(} \beta(y) - \sum_{j <  i} \psi_j(y_j)  - \sum_{j >  i} g_j(y_j) \big{)}, 
\]
we find that $\psi_i$ are convex, and 
\begin{align*}
\beta(y) \le \sum_{i=1}^d \psi_i (y_i) \quad \forall y \in \R^d, \quad\text{and}\quad\beta(y) = \sum_{i=1}^d \psi_i (y_i) \quad \pi^2 - a.e. \, y,
\end{align*}
where  $\pi^2= p^y_\#\pi$ as usual. Let us define the contact set
\begin{align*}
G= \{(x,y) \in \R^{2d} \,|\, f^\oplus(x) + h(x) \cdot (y-x) - c(x,y) =\beta(y)= \psi^\oplus(y)\}.
\end{align*}
Let us slightly refine $G$ as follows, and define its projections:
\begin{align*}
H&= \{(x,y) \in G \,|\, \psi_i \text{ is differentiable at } y_i \text{ for every } i \in S\}\\
X_H &= \{x \in \R^d \,|\,\exists y\,\, s.t.\, (x,y) \in H\},  \ Y_H = \{y \in \R^d \,|\,\exists x\,\, s.t.\, (x,y) \in H\}.
\end{align*}
As $\nu_i \ll \mathcal{L}$ for all $i \in S$ and $\psi_i$ is convex, $\pi(H)=1$ for any {\sf VMOT} $\pi$. 

Now we claim that for any $(x,y) \in H$ and $i \in S$,
\begin{align}\label{diffy}
&\frac{\partial \psi_i}{\partial y_i}(y_i), \ \frac{\partial \beta}{\partial y_i}(y) \ \text{and} \ \frac{\partial c}{\partial y_i}(x,y) \ \text{exist, and} \\
\label{diffyeq}& \frac{\partial \psi_i}{\partial y_i}(y_i)= \frac{\partial \beta}{\partial y_i}(y)= h_i(x) - \frac{\partial c}{\partial y_i}(x,y).  
\end{align}
Of course $\frac{\partial \psi_i}{\partial y_i}(y_i)$ exists by definition of $H$. Also from the definition 
\[
\psi_i(y_i) = \sup_{y_j , j \neq i} \big{(} \beta(y) - \sum_{j \neq i} \psi_j (y_j) \big{)}
\]
we see that the first equality in \eqref{diffyeq} holds since the supremum is attained at $y=(y_i)_i$ and both $\psi_i$ and $\beta$ are convex. Then in turn since the supremum in \eqref{1beta} is attained at $(x,y)$ and both $\beta$ and $y \mapsto h(x) \cdot (y-x) - c(x,  y)$ are convex, the second equality holds as well.

For each $x \in X_H$ define the slice set $H_x := \{y \,|\, (x,y) \in H\}$, and for $y=(y_1,...,y_d) \in \R^d$ define the projection $y_S \in \R^{|S|}$ of $y$ to be the collection of those $y_i$'s with $i \in S$. Now by \eqref{diffyeq} and the twist condition \eqref{twist}, we see that if $y^1, y^2 \in H_x$ and $y^1_S = y^2_S$, then we must have $y^1=y^2$. This implies that we can define a function $F_x : \R^{|S|} \to \R^{d-|S|}$ such that the set $H_x$ is contained in the graph of $F_x$ for every $x \in X_H$. Finally as $\pi(H) =1$, we have $\pi_x(H_x)=1$ for $\pi^1$-a.e.\,$x$, completing the proof.
\end{proof}

\begin{remark}
In Theorem \ref{mongesolution} if one could obtain a dual maximizer $(f_i, -g_i,h_i)_i$ where all $g_i$'s and  $y \mapsto c(x,y)$ are differentiable, then the family of maps $\{F_x\}_x$ could be directly obtained by the above differential identity in $y$; see the next example. But instead the semiconcavity was assumed in the theorem in order to deal with more general costs.
\end{remark}

\begin{example}
Set $d=2$ and $c(x,y) = - y_1y_2$, so that we consider the covariance maximization problem of $Y_1, Y_2$. Let $(f_i, -g_i,h_i)$ be a {\sf PDM}. In this case we can assume that $g_1,g_2$ are convex, since we can replace $g_1$ by the following Legendre transform $\tilde g_1$ (and similarly for $g_2$)
\[
\tilde g_1 (y_1) := \sup_{x_1,x_2,y_2}\bigg{(}\sum^2_{i=1}\big{(} f_i(x_i) + h_i(x_1,x_2) \cdot (y_i-x_i)\big{)}  - c(x,y) -g_2(y_2) \bigg{)}.
\]
So $g_1,g_2$ are differentiable a.e., and the differential identity in $y_1$ reads
\[ y_2 =   g_1 '(y_1) - h_1(x_1,x_2) 
\]
which represents the function $y_2=F_x(y_1)$ in Theorem \ref{mongesolution}. Interestingly, even though $c(x,y)$ has no dependence on $x$, $F_x$ still depends on $x$ via the ``trading strategy" $h_1$. This is due to the randomness of $X_1,X_2$ and the martingale constraint $\E[Y|X]=X$. Note that if $X_1, X_2$ are nonrandom, then one can take $h_1 \equiv 0$ and the function becomes the monotone map. Nonetheless, observe that even if $X_1, X_2$ are random, the dependence of the functions $y_2=F_x(y_1)$ on $x$ is mild; for any given marginal $(\vec\mu, \vec\nu)$, their graphs are merely translations in $y_2$ direction. 
\end{example}

Theorem \ref{extreme} and Theorem \ref{mongesolution} discussed the optimal structure of $(\pi_x)_x$ for {\sf VMOT} $\pi= \pi_x \otimes \pi^1$. Would $\pi^1= p^x_\# \pi$ and $\pi^2=p^y_\# \pi$ also have some optimality property? To discuss, recall \eqref{alpha1}, \eqref{H} where we defined $y \mapsto H(x,y)$ to be the convex envelope of $y \mapsto c(x,y) + g^\oplus (y)$ and set $\alpha(x) = H(x,x)$. Also recall $\beta$ is defined by \eqref{1beta}. The following theorem shows the optimality of $\pi^1, \pi^2$ with respect to the costs $\alpha, \beta$ respectively from a classical duality theorem of Kellerer \cite{Ke84}.
\begin{theorem}\label{otduality}
Let $\pi \in {\rm VMT}(\vec\mu, \vec\nu)$ be a minimizer for \eqref{VMOT} and let $(f_i, -g_i, h_i)_{i \in [d]}$ be a dual maximizer. Then $\pi^1, \pi^2$  satisfy the following:  \begin{align}
\label{otdual1}
\sum_i f_{i}(x_i) \le \alpha(x)  \ \ \forall x \in \R^d,   &\text{ and }   \sum_i f_{i}(x_i) = \alpha(x) \ \ \pi^1 - a.e. \, x,  \\
\label{otdual2}
\sum_i g_{i}(y_i) \ge \beta(y)   \ \ \forall y \in \R^d,  &\text{ and }   \sum_i g_{i}(y_i) = \beta(y) \ \  \pi^2 - a.e. \,y.
\end{align}
\end{theorem}
\begin{proof} Recall $f^\oplus(x) + h(x) \cdot (y-x)  \le c(x,y) + g^\oplus(y)$. Then by convexity,
 \begin{align}\label{usual}
 f^\oplus(x) + h(x) \cdot (y-x) \le  H(x,y) \le c(x,y) + g^\oplus(y).
\end{align}
By setting $y=x$ in the first inequality, we get $ f^\oplus(x) \le H(x,x) =:\al (x)$, and notice from the first inequality that if $ f^\oplus(x)=\al (x)$ then $h(x)$ must belong to the subdifferential of the convex function $y \mapsto H(x,y)$ at $x$. Now since \eqref{usual} must hold as equalities for $\pi$-a.e. $(x,y)$ and $\delta_x \le_c \pi_x$, we deduce that $ f^\oplus(x)=\al (x)$ must hold $\pi^1$-a.e. $x$. This yields \eqref{otdual1}.

The proof of \eqref{otdual2} is similar. By definition of $\bt$, we have
 \begin{align}\label{usual1}
 f^\oplus(x) + h(x) \cdot (y-x)- c(x,y) \le \bt(y) \le g^\oplus(y)
\end{align}
and moreover \eqref{usual1} must hold as equalities $\pi$-a.s.. This yields \eqref{otdual2}.
\end{proof}

\begin{example}\label{mmotisot}
In Example \ref{ex1} and \ref{ex2}, the optimality of $\pi$ did not imply any constrained structure on $\pi^1$, since $\pi^1$ could be any coupling of $(\mu_1,...,\mu_d)$ there. In this example we shall see that for some data (i.e. cost and marginals), $\pi^1$ may have to be constrained. As in the previous example, let $d=2$ and $c(x,y) = - y_1y_2$.

It is well known that the monotone coupling, which is the coupling  concentrated on an increasing graph, minimizes the cost (i.e. maximizes $\E[Y_1Y_2]$) among all couplings of $\nu_1, \nu_2$. But in our {\rm VMOT} setting, the monotone coupling may not be feasible to be $\pi^2$ since we have to satisfy $\pi^1 \le_c \pi^2$ but there may not exist such a $\pi^1$. Still,  $\pi^2$ may try to be as close as possible to the monotone coupling, and this in turn may affect the structure of $\pi^1$. We explore this phenomenon with explicit examples.

Firstly, as a toy model let the marginals be those in Example \ref{ex2}: $\mu_1 = \mu_2 = {\rm Leb}\big{|}_{[-1/2, 1/2]}$, $\nu_1= \nu_2 = \frac{1}{2}{\rm Leb}\big{|}_{[-1, 1]}$ and let $\pi$ be a {\rm VMOT}. As maximizing $\E[Y_1Y_2]$ is equivalent to minimizing $\E\frac{1}{2}|Y_1-Y_2|^2$, it is better for $\pi^2$ to stay near the diagonal $\Delta := \{(y_1,y_2) \,|\, y_1=y_2 \}$. With this simple data $(\vec\mu, \vec\nu)$  it is feasible for $\pi^2$ to be supported on $\Delta$, and since $\pi^1 \le_c \pi^2$, this implies $\pi^1$ must also be supported on $\Delta$. This uniquely determines $\pi^1, \pi^2$, and any martingale measure $\pi$ connecting $\pi^1,\pi^2$ is optimal. Thus a {\rm VMOT} $\pi$ is nonunique but $\pi^1, \pi^2$ are unique.

For another example, let $\mu_1, \mu_2$ be arbitrary but fixed. The main goal is to build a {\rm VMOT} $\pi^*$ such that its first marginal $\pi^{*1}$ is the monotone coupling of arbitrarily given marginals $\mu_1, \mu_2$. To this end, set $g_1(y_1) = \frac{1}{2}|y_1|^2$, $g_2(y_2) = \frac{1}{2}|y_2|^2$, so that $c(x,y) + g_1(y_1) + g_2(y_2) = \frac{1}{2}|y_1-y_2|^2$. Now since $\frac{1}{2}|y_1-y_2|^2$ is already convex in $y$, the first derived cost $\alpha$ appearing in Theorem \ref{otduality} must be the same, i.e. $\alpha(x) =  \frac{1}{2}|x_1-x_2|^2$. This implies $h(x) = \nabla \alpha(x) = (x_1 - x_2, x_2 - x_1)$, thus $y \mapsto \alpha(x) + h(x) \cdot (y-x)$ is the tangent plane to $y \mapsto c(x,y) + g_1(y_1) + g_2(y_2)$ at $x$. In this simple case,  it is easy to see what the contact set is given by:
\begin{align*}
G:= &\{ (x,y) \,|\, \alpha(x) + h(x)\cdot(y-x) = c(x,y) + g_1(y_1) + g_2(y_2)\}\\
=&\{ (x,y) \,|\, y-x = \lambda(1,1) \text{ for some $\lambda \in \R$}    \}.
\end{align*} 
Now choose a martingale measure $\pi^*$ which satisfies that $\pi^*(G)=1$, for which its first induced coupling $\pi^{*1} = p^x_\# \pi^*$ is the monotone coupling of $\mu_1, \mu_2$. This implies there exist Kantorovich potentials $f_1, f_2$ such that
\begin{align}
\label{brenierineq}&f_1(x_1)+ f_2(x_2)  \le \alpha(x)  \quad  \forall x=(x_1, x_2),\ \text{and}\\
\label{breniereq}&f_1(x_1)+ f_2(x_2)  = \alpha(x)  \quad  \pi^{*1} - a.e.\, x.
\end{align}
Let $\nu_1^*, \nu_2^*$ be the one-dimensional marginals of $\pi^{*2}= p^y_\# \pi^*$. Then the fact $\pi^*(G)=1$ and \eqref{brenierineq}, \eqref{breniereq} imply that $\pi^*$ is optimal in the class ${\rm VMT}(\mu_1,\mu_2,\nu_1^*,\nu_2^*)$. Moreover, any optimal $\pi \in {\rm VMT}(\mu_1,\mu_2,\nu_1^*,\nu_2^*)$ must satisfy $\pi^1 = \pi^{*1}$ since \eqref{brenierineq},  \eqref{breniereq} must hold for $\pi^1$ as well.
\end{example}
In this example the choices $g_i(y_i) = \frac{1}{2}|y_i|^2$ were made with no specific reason, and due to this the $\pi^*$ constructed there may not fit to a prescribed marginal data $(\vec \mu, \vec \nu)$. Theorem \ref{main2} asserts that for any cost and irreducible marginal data one can obtain a dual optimizer such that the contact set \eqref{ptwiseeq} can accommodate all {\sf VMOT}s. Moreover, we have illustrated that some careful analysis of the contact set can provide information on the  structure of {\sf VMOT}s, as shown in Theorem \ref{extreme}, \ref{mongesolution}.

Lastly, motivated by this example, we address the following question: what conditions on the data $(\vec \mu, \vec \nu)$ and cost $c(x,y)$ shall impose on the induced optimal couplings $\pi^1$, $\pi^2$ to have some specific structures, e.g.  $\pi^1$ or $\pi^2$ has to lie on  ``small" sets? We leave this to future research.

\appendix

\section{Compactness of the convex potentials}\label{convexpotential}
We prove Proposition \ref{conv} which we restate for reader's convenience.
\begin{proposition}\label{convv} Let $(\mu_i,\nu_i)_{i \in [d]}$ be  irreducible pairs of probabilities with domains $(I_i,J_i)_i$. Let $a \in I$, $C \in \R$. Consider the following class of functions $\Lambda=\Lambda(a,C,\vec\mu,\vec\nu) $ where every $\chi \in \Lambda$ satisfies the following:
\begin{enumerate}
\item $\chi$ is a real-valued convex function on $J$,
\item $\chi \ge 0$ and $\chi(a)=0$,
\item $\int \chi \,d(\nu^\otimes - \mu^\otimes) \le C$.
\end{enumerate}
Then $\Lambda$ is locally bounded in the following sense: for each compact subset $K$ of $J$, there exists $M=M(K)$ such that $\chi \le M$ on $K$ for every $\chi \in \Lambda$. Moreover, for any sequence $\{\chi_n\}_n$ in $\Lambda$ there exists a subsequence $\{ \chi_{n_j} \}_j$ of $\{ \chi_n\}_n$ and a real-valued convex function $\chi$ on $J$ such that $\lim_{j \to \infty} \chi_{n_j} (x) = \chi(x)$ for every $x \in J$, and  the convergence is uniform on every compact subset of $I$.
\end{proposition}
The following class of two-way martingales will also be used often.
\begin{definition}
(1) We define $\zeta = \zeta_{x \to (y^-, y^+)}$ to be a simple martingale measure in $\R \times \R$ such that $p^x_\# \zeta$ is supported at a point (say $x$), and $p^y_\#\zeta$ is supported at two distinct points (say $y^-, y^+$). Then necessarily $y^- < x < y^+$, and $\zeta$ can be written as
 \[\zeta_{x \to (y^-, y^+)} = \frac{y^+ - x}{y^+ - y^-} \delta_{(x, y^-)} + \frac{x - y^- }{y^+ - y^-} \delta_{(x, y^+)}.
 \]
Then we define $\mathfrak{m}$ to be the set of all such simple martingale measures.
\\
(2) Let $\mathfrak{M}$ be the set of martingale measures in $\R \times \R$ such that any $\pi =\pi_x \otimes \pi^1 \in \mathfrak{M}$ satisfies $\pi_x \otimes \delta_x \in \mathfrak{m}$ for $\pi^1$--a.e. $x$, where $\pi^1=p^x_\#\pi$.
\smallskip

\noindent (3) We define $\mathfrak{M}(\mu,\nu) := \mathfrak{M} \cap {\rm MT}(\mu,\nu)$.
\end{definition}
\begin{proof}[Proof of Proposition \ref{convv}]
{\bf Step 1.}  In this step we will prove for the simplest marginal case where $d=1$ and there exists $\zeta = \zeta_{x \to (y^-, y^+)}    \in \mathfrak{m}$ such that $\mu = p^x_\# \zeta$ and  $\nu = p^y_\# \zeta$, that is $\mu = \delta_x$ and $\nu=  \frac{y^+ - x}{y^+ - y^-} \delta_{y^-} + \frac{x - y^- }{y^+ - y^-} \delta_{y^+}$. For $x,y \in \R$, $x \neq y$, define the open / closed transport rays
\begin{align*}
\rrbracket x,y  \llbracket = \{ (1-t)x+ty \,|\, 0<t<1\}, \ \llbracket x,y  \rrbracket = \{ (1-t)x+ty \,|\, 0\le t \le 1\}.
\end{align*}

Suppose there is $a \in I=\, \rrbracket y^-,y^+ \llbracket$ and $C>0$ such that $\chi(a)=0$, $\chi \ge 0$  and $\int \chi \,d(\nu-\mu) \le C$ for every $\chi \in \Lambda$. First, consider the case $a = x$. Then it is trivial that there is a constant $M>0$ such that $\chi(y^-) \le M$ and $\chi(y^+) \le M$ for all $\chi$, since $\int \chi \,d(\nu-\mu) = \int \chi \,d\nu = \frac{y^+ - x}{y^+ - y^-} \chi({y^-}) + \frac{x - y^- }{y^+ - y^-} \chi({y^+}) \le C$. By convexity of $\chi$, we find $\chi \le M$ on $\llbracket y^-,y^+  \rrbracket$ for every $\chi \in \Lambda$, proving the theorem in this most simple case.

Next, consider the case that $a$ is away from $x$, e.g. $y^- < a < x$. The idea is to find probability measures $\tilde \mu, \tilde \nu$ where $\mu \le_c \tilde \mu \le_c \tilde \nu \le_c \nu$, and to find a small but positive constant $\theta$ (which does not depend on $\chi \in \Lambda$) such that the martingale measure $\eta := \theta\, \zeta_{a \to (y^-, y^+)}$ satisfies $\eta \le \tilde \pi$ for some $\tilde \pi \in {\rm MT}(\tilde \mu, \tilde \nu)$, i.e., for some martingale measure $\tilde \pi$ having marginals $\tilde \mu, \tilde \nu$. If we can do so, then since $\tilde \pi - \eta$ is still a martingale, we infer from $\int \chi \,d(\tilde \nu-\tilde \mu) \le C$ that we have $\int \chi \,d(p^y_\# \eta-p^x_\# \eta) \le C$. Then as $\eta := \theta\, \zeta_{a \to (y^-, y^+)}$, we are reduced to the previous case and the proposition follows. Now for this, we  can simply take $\tilde \nu = \nu$ and $\tilde \mu = p^y_\# \zeta_{x \to (a, y^+)}$. Since $y^- < a < x$, clearly $\tilde \mu \le_c \nu$. Then notice we can take $\theta = \tilde \mu (a)= \frac{y^+ - x}{y^+ - a} $, implying $\tilde \pi = \eta + \tilde \mu (y^+) \delta_{(y^+, y^+)}$; see figure \eqref{step1martingale}.

We can describe this procedure in terms of Brownian motion and stopping time, namely, let $(B^x_t)_t$ denote a Brownian motion starting at $x$, and define the stopping time $\tau = \inf \{t \ge 0 \,|\, B^x_t \notin ]y^-, y^+[ \}$. Then ${\rm Law}(B^x_{\tau}) = \nu$. Next, define  $\tau^1 = \inf \{t \ge 0 \,|\, B^x_t \notin ]a , y^+[ \}$ and $\tilde \mu = {\rm Law}(B^x_{\tau^1})$. Then as $\tau^1 \le \tau$ we have $\tilde \mu \le_c \nu$. Finally, define $\tau^2$ as the following: if $B^x_{\tau^1} = y^+$ then let $\tau^2 = \tau^1$, and if $B^x_{\tau^1} = a$ then let $\tau^2 = \tau ^1 + \inf \{t \ge 0 \,|\, B^a_t \notin ]y^- , y^+[ \}$. We may define $\tilde \nu = {\rm Law}(B^x_{\tau^2})$, but in this case $\tilde \nu = \nu$ simply because $\tau^2 = \tau$. 

The case $x < a < y^+$ can be treated similarly and  Step 1 is done.

\begin{figure}
    \begin{subfigure}[b]{0.48\columnwidth}
        \includegraphics[width=\textwidth]{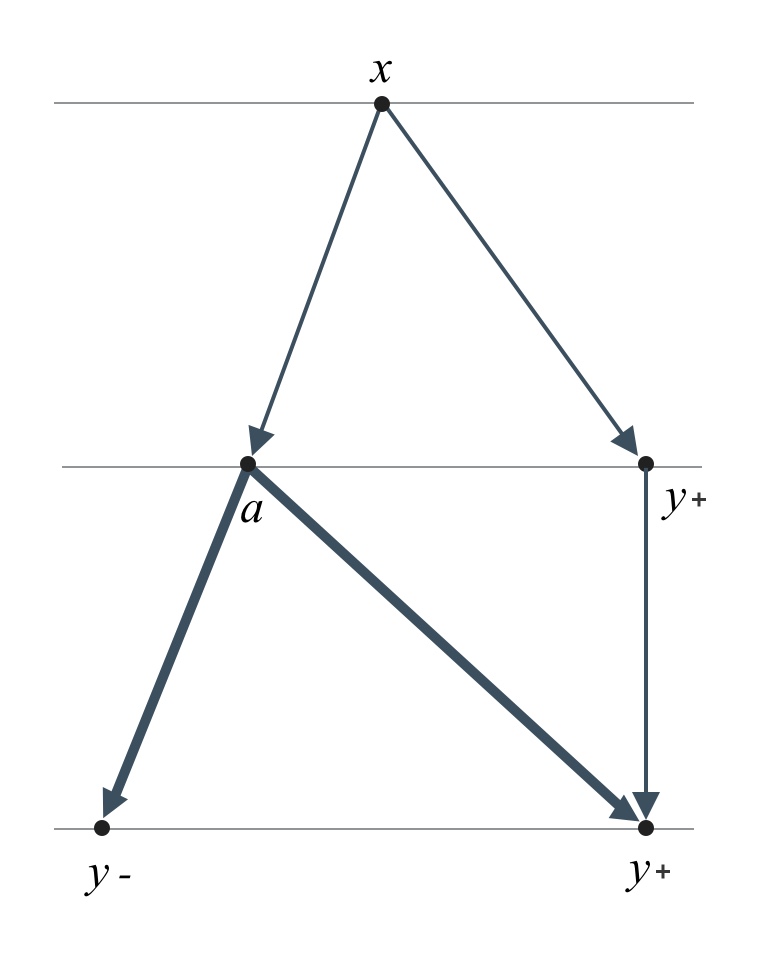}
      \caption{Step 1}
        \label{step1martingale}
    \end{subfigure}
    \hfill
    \begin{subfigure}[b]{0.48\columnwidth}
        \includegraphics[width=\textwidth]{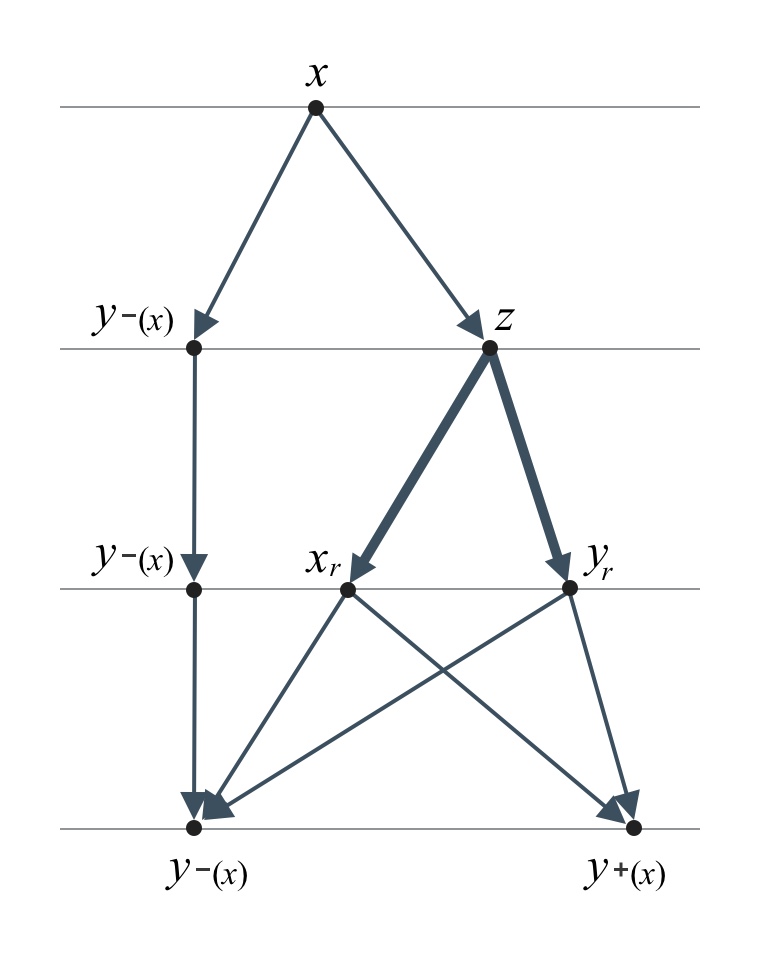}
      \caption{Step 3}
      \label{step3martingale}
    \end{subfigure}
    \caption{Martingale constructions in Step 1 and 3}
  \label{step13}
\end{figure}

{\bf Step 2.} In this step we will prove for the case where $(\mu,\nu)$ are the marginals of a martingale measure of the form $\pi := \sum_{i=1}^N p_i\zeta_{x_i \to (y^-_i, y^+_i)}$. Here $N \ge 2$, $p_i >0, \sum p_i =1$, and the family $\zeta_{x_i \to (y^-_i, y^+_i)} \in \mathfrak{m}$, $i=1,2,...,N$ is assumed to  satisfy the following {\em chain condition}: if we let $l_k = \min_{i \le k} y^-_i$, $r_k = \max_{i \le k} y^+_i$, then
\begin{align*}
\rrbracket l_k,r_k  \llbracket \,\, \cap \,\, \rrbracket y^-_{k+1}, y^+_{k+1}  \llbracket \,\, \neq \,\, \emptyset, \quad \forall k=1,...,N-1.
\end{align*}
The chain condition immediately implies that $(\mu,\nu)$ is irreducible. As it is enough to prove this step with $N=2$ let us assume $N=2$. Without loss of generality assume $a \in \rrbracket y^-_1,y^+_1 \llbracket$ such that $\chi(a)=0$, $\chi \ge 0$  and $\int \chi \,d(\nu-\mu) \le C$ for every $\chi \in \Lambda$. Then since $p_1 \zeta_{x_1 \to (y^-_1, y^+_1)} \le \pi$, by Step 1 we see that there exist $M_1, M'_1 >0$ 
such that for every $\chi \in \Lambda$
\begin{enumerate}
\item $ \chi \le M_1$ on  $\llbracket y^-_1, y^+_1 \rrbracket$, and
\item $\nabla \chi \le M'_1$ on  $\llbracket y^-_1 + \epsilon , y^+_1 -\epsilon \rrbracket$, for any fixed $\epsilon>0$.

\end{enumerate}
$M'_1$ may depend on $\epsilon$ but notice it does not depend on $\chi  \in \Lambda$. Note that (1) follows from Step 1 and (2) follows from (1) and the fact that $\chi$ is convex. Now take $\epsilon$ small enough such that 
\begin{align*}
\rrbracket y^-_1 + \epsilon , y^+_1 -\epsilon \llbracket \,\, \cap \,\, \rrbracket y^-_2 + \epsilon , y^+_2 -\epsilon  \llbracket \,\, \neq \,\, \emptyset.
\end{align*}
Fix a point $b \in \,\, \rrbracket y^-_1 + \epsilon , y^+_1 -\epsilon \llbracket \,\, \cap \,\, \rrbracket y^-_2 + \epsilon , y^+_2 -\epsilon  \llbracket$ and take any $\chi \in \Lambda$. Let $L_\chi$ be an affine function which supports the convex function $\chi$ at $b$, that is $L_\chi(b) = \chi(b), \nabla L_\chi(b) \in \partial \chi(b)$. Note that, by (1) and (2), $L_\chi(b)$ and $\nabla L_\chi(b)$ are bounded independently of $\chi$. Let $\tilde \chi = \chi - L_\chi$. Then $\tilde \chi(b) =0$, $\tilde \chi \ge 0$  and $\int  \tilde  \chi \,d(\nu-\mu) \le C$, hence by the observation (1), (2) and the Step 1 we deduce there exists $M_2 > 0$ 
such that for every $\chi \in \Lambda$ we have $ \chi \le M_2$ on  $\llbracket y^-_2, y^+_2 \rrbracket$. Generalization to arbitrary $N \in \N$ is immediate through the chain condition, and this completes Step 2.\\

{\bf Step 3.} In this step we will prove for a general irreducible pair $(\mu,\nu)$ with domain $(I,J)$ but still in dimension $d=1$. In this case, it is well-known that there exists a probability measure $\tilde \mu$ such that $\tilde \mu$ is absolutely continuous with respect to Lebesgue measure, $\mu \le_c \tilde \mu \le_c \nu$, and $(\tilde \mu,\nu)$ is irreducible with the same domain $(I,J)$.  One way to see this is the following: consider $u_\mu, u_\nu$, the potential functions of $\mu, \nu$. They coincide outside of $I$ while  $u_\mu < u_\nu$ in $I$ by irreducibility. Select a convex function $u$ such that $u_\mu = u = u_\nu$ outside of $I$ while $u_\mu <u < u_\nu$ in $I$. Then $\tilde \mu$ is taken to be the second derivative measure of $u$,  and by selecting sufficiently smooth $u$ one gets the desired  absolute continuity.

Hence from now on we will assume without loss of generality that $\mu \ll \mathcal{L}^1$. Also note that we can assume $\mu \wedge \nu =0$  since  the pair $(\mu - \mu \wedge \nu, \nu - \mu \wedge \nu)$ is irreducible with domain $(I,J)$ and $\int \chi \, d \big{(}(\nu - \mu \wedge \nu) - (\mu - \mu \wedge \nu)\big{)} \le C$. Thus, we will assume $\mu \ll \mathcal{L}^1$ and $\mu \wedge \nu =0$.

Our strategy is to reduce this situation to Step 1 and Step 2. For this, we begin by choosing  a martingale measure $\pi \in \mathfrak{M}(\mu,\nu)$. We can obtain such a $\pi$ by solving the {\sf MOT} problem w.r.t. the cost $c(x,y) = \pm |x-y|$\cite{HoKl12}, \cite{HoNe11}, \cite{bj}. Let $x \mapsto y^- (x)$, $x \mapsto y^+ (x)$ be two functions supporting $\pi = \pi_x \otimes \mu$, that is  $\pi_x = \frac{y^+(x) - x}{y^+(x) - y^-(x)} \delta_{y^-(x)} + \frac{x - y^-(x) }{y^+(x) - y^-(x)} \delta_{y^+(x)}$ $\mu$ - a.e. \,$x$.
 
Now set $\Gamma := \supp \pi$ such that $\pi(\Gamma)=1$ and for any $(x,y) \in \Gamma$ and $r >0$, the restricted measure $ \pi \big{|}_{B(x,r) \times B(y, r)}$ has strictly positive mass by the property of support. Pick $(x_0,y_0) \in \Gamma$ and  suppose $x_0 <y_0$. Fix $z \in \rrbracket x_0, y_0 \llbracket$ and let  $r >0$  be so small that  $x_r := x_0 + r <z< y_r := y_0 - r$. Now we claim the following: there exist positive measures $\tilde \mu, \tilde \nu, \tilde \mu_A, \tilde \nu_A$ satisfying $\tilde \mu_A \le_c \tilde \mu \le_c \tilde \nu \le_c \tilde \nu_A$, martingale measures $\tilde \pi \in {\rm MT}(\tilde \mu, \tilde \nu)$, $\pi_A \in {\rm MT}(\tilde \mu_A, \tilde \nu_A)$ such that $\pi_A \le \pi$, and a  positive constant $\tilde \theta$ 
such that the measure $\eta :=\tilde  \theta\, \zeta_{z \to (x_r, y_r)}$ satisfies $\eta \le \tilde \pi$. Notice this yields $\int \chi\, d(\nu - \mu) \ge \int \chi\, d(p^y_\# \eta - p^x_\# \eta)$ as desired. See figure \eqref{step3martingale}.
 
To see this, note that the set $A:= \{x \in B(x_0,r)  \,|\, y^+(x) \in B(y_0,r) \}$ has $\mu$ -- positive measure since $ \pi  \big{(} B(x_0,r) \times B(y_0, r)\big{)}>0$. For each $x\in A$, consider the martingale measure $\zeta_{x \to (y^-(x), y^+(x))}$. Recall $y^-(x) < x < x_r < z < y_r< y^+(x)$.  Now we apply the same idea as in Step 1: Let $\tau_x = \inf \{t \ge 0 \,|\, B^x_t \notin ]y^-(x), y^+(x)[ \}$ and $\tau^1_x = \inf \{t \ge 0 \,|\, B^x_t \notin ]y^-(x), z[ \}$. Then define $\tau^2_x$ as follows: if $B^x_{\tau^1_x} = y^-(x)$ then $\tau^2_x = \tau^1_x$, and if $B^x_{\tau^1_x} = z$ then $\tau^2_x = \tau ^1_x + \inf \{t \ge 0 \,|\, B^z_t \notin ]x_r, y_r[ \}$. Then clearly 
\begin{align}\label{convexorder1}
\delta_x \le_c {\rm Law}(B^x_{\tau^1_x}) \le_c {\rm Law}(B^x_{\tau^2_x}) \le_c {\rm Law}(B^x_{\tau_x}).
\end{align}
With $\theta = \theta(x,z) := {\rm Law}(B^x_{\tau^1_x}) (\{z\}) = \frac{x-y^-(x)}{z-y^-(x)}$,  $\theta \, \zeta_{z \to (x_r, y_r)}$ satisfies
\begin{align}\label{convexorder2}
\theta \, \zeta_{z \to (x_r, y_r)} \le {\rm Law}(B^x_{\tau^1_x}, B^x_{\tau^2_x}),
\end{align}
where ${\rm Law}(B^x_{\tau^1_x}, B^x_{\tau^2_x})$ is the joint law of $B^x_{\tau^1_x}$ and $B^x_{\tau^2_x}$ hence is a martingale measure with marginals ${\rm Law}(B^x_{\tau^1_x})$, ${\rm Law}(B^x_{\tau^2_x})$. Finally, set
 \begin{align*}
 \tilde \theta (z) = \int_A \theta(x,z) \,d\mu(x).
  \end{align*}
Notice $ \tilde \theta (z) > 0$.  We claim that the measure $\tilde \theta \, \zeta_{z \to (x_r, y_r)}$ satisfies the desired property. To see this, consider the restricted measure $\pi_A :=  \pi  \big{|}_{ A \times \R} \in {\rm MT}(\tilde \mu_A, \tilde \nu_A)$ where $\tilde \mu_A := p^x_\# \pi_A$, $\tilde \nu_A := p^y_\# \pi_A$. Define $\tilde \mu = \int_A {\rm Law}(B^x_{\tau^1_x}) \,d\mu(x)$, $\tilde \nu = \int_A {\rm Law}(B^x_{\tau^2_x}) \,d\mu(x)$. Then in view of \eqref{convexorder1} we have $\tilde \mu_A \le_c \tilde \mu \le_c \tilde \nu \le_c \tilde \nu_A$, and by \eqref{convexorder2} and  $\eta :=\tilde  \theta\, \zeta_{z \to (x_r, y_r)}$ we have $\eta \le \tilde \pi$ for some $\tilde \pi \in {\rm MT}(\tilde \mu, \tilde \nu)$, where $\tilde \pi$ is driven by the Brownian motion $B^x$ between times $\tau^1_x$ and $\tau^2_x$ for $\mu\big{|}_A$-- a.e. $x$.

In summary, we have shown that for any $(x,y) \in \Gamma$ and $z \in \rrbracket x, y \llbracket$ and any sufficiently small $r>0$, we can find $\tilde \mu, \tilde \nu, \tilde \mu_A, \tilde \nu_A$ where $\tilde \mu_A \le_c \tilde \mu \le_c \tilde \nu \le_c \tilde \nu_A$ and martingale measures $\tilde \pi \in {\rm MT}(\tilde \mu, \tilde \nu)$, $\pi_A \in {\rm MT}(\tilde \mu_A, \tilde \nu_A)$ such that $\pi_A \le \pi$, and a  positive constant $\tilde \theta$ 
such that the measure $\eta :=\tilde  \theta\, \zeta_{z \to (x_r, y_r)}$ satisfies $\eta \le \tilde \pi$. Now we make the following observation. Since $\pi(\Gamma)=1$,  for any $z \in I$ there exists $(x_z, y_z) \in \Gamma$ such that $z \in \rrbracket x_z, y_z \llbracket$, since otherwise we would have $u_\mu (z) = u_\nu (z)$, contradicting to the irreducibility of $(\mu,\nu)$. We choose such $(x_z, y_z) \in \Gamma$ for each $z \in I$ so that $\{ \,\rrbracket x_z, y_z \llbracket \, \}_{z \in I}$ is an open cover of $I$. Hence we can find a sequence $(z_n)_{n \in \N}$ such that $\{ \,\rrbracket x_{z_n}, y_{z_n} \llbracket \, \}_{n \in \N}$ is an open cover of $I$, and if $l_k := \min_{n \le k} x_{z_n}$, $r_k := \max_{n \le k} y_{z_n}$, then $\rrbracket l_k,r_k  \llbracket \,\, \cap \,\, \rrbracket x_{z_{k+1}}, y_{z_{k+1}}  \llbracket \,\, \neq \,\, \emptyset$ for all $ k\in \N$, namely, the chain condition.

Now Step 2 implies that for each compact interval $K \subset I$, there exists $M_K>0$ such that $\chi \le M_K$ on $K$ for every $\chi \in \Lambda$, which is the content of the proposition when $\nu$ has no mass on the boundary of $I$. To see this,  observe that for a compact $K$ there is $N \in \N$ such that $\{ \,\rrbracket x_{z_n}, y_{z_n} \llbracket \, \}_{n \le N}$ is an open cover of $K$, so for small $r>0$ the shrunken covering $\{ \,\rrbracket x_{z_n} + r, y_{z_n} - r \llbracket \, \}_{n \le N}$ still covers $K$ and satisfies the chain condition, hence Step 2 applies. This completes Step 3 for the case $\nu(I)=1$. Now the remaining case when $\nu$ assigns positive mass on the boundary of $I$ is simple. Write $I= ]a,b[$ and suppose $b \in \R$ and $\nu(b)>0$. Recall that $\pi \in \mathfrak{M}(\mu,\nu)$ and $x \mapsto (y^- (x), y^+ (x))$ is the graph of $\pi$. Let $A := \{x \,|\, y^+(x) =b \}$. $\nu(b)>0$ implies $\mu(A) >0$. Take any $x \in A$ satisfying  $(x,b) \in \Gamma$. Then the previous argument shows that for any $r>0$ and $z \in ]x+r , b[$, we can find measures $\tilde \mu, \tilde \nu, \tilde \mu_A, \tilde \nu_A$ where $\tilde \mu_A \le_c \tilde \mu \le_c \tilde \nu \le_c \tilde \nu_A$, martingale measures $\tilde \pi \in {\rm MT}(\tilde \mu, \tilde \nu)$, $\pi_A \in {\rm MT}(\tilde \mu_A, \tilde \nu_A)$ such that $\pi_A \le \pi$, and $\tilde \theta > 0$ 
such that $\eta :=\tilde  \theta\,\zeta_{z \to (x+r, b)}$  satisfies $\eta \le \tilde \pi$. Step 2 then yields for any $t \in I$, there exists $M >0$ such that $\chi \le M$ on $[t,b]$ for every $\chi$.

The case $\nu(a) >0$ can be treated similarly, and finally note that if both $\nu(a) >0$ and $\nu(b) >0$ then we can conclude that 
there exists $M >0$ such that $\chi \le M$ on $[a,b]$ for every $\chi \in \Lambda$, because $[a,b]$ can be covered with finitely many appropriate intervals that we have discussed so far. The proof is therefore complete for the one-dimension case.\\

{\bf Step 4.} In this step we will prove the local uniform bound $\chi \le M$ for general $d \in \N$. We will see how the ideas from the previous steps can be applied in the same way, but with more notational difficulty; this is why we take a different approach than \cite{bnt} in establishing the local uniform bound of $\chi$ for $d=1$ in the previous steps.

Fix $d \in \N$. As in Step 1, we begin with the most simple case where for each $i \in [d]$ there is $\zeta_i = \zeta_{x_i \to (y^-_i, y^+_i)}    \in \mathfrak{m}$ such that $\mu_i = p^x_\# \zeta_i$, $\nu_i = p^y_\# \zeta_i$, that is, $\mu_i = \delta_{x_i}$ and $\nu_i=  \frac{y^+_i - x_i}{y^+_i - y^-_i} \delta_{y^-_i} + \frac{x_i - y^-_i }{y^+_i - y^-_i} \delta_{y^+_i}$. Let  $\mu := \otimes_i \mu_i$, $\nu := \otimes_i \nu_i$ and $\zeta := \otimes_i \zeta_i$. Note that the martingale measure $\zeta \in \cP(\R^{2d})$ has $\mu, \nu$ as its $d$-dimensional marginals. Now we want to prove that there exists $M >0$ such that $\chi \le M$ on $J$ for all $\chi \in \Lambda$, where $J= {\rm conv}(\supp \nu)$ is the convex hull of the support of $\nu$. Notice that $J$ is a closed rectangle  in $\R^d$ and $\supp \nu$ is the set of its vertices, consisting of $2^d$ points.

Suppose there is $a \in I$ (where $I = {\rm int} \,J$) and $C>0$ such that $\chi(a)=0$, $\chi \ge 0$  and $\int \chi \,d(\nu-\mu) \le C$ for every $\chi \in \Lambda$. First, consider the case $a = x$, where $x = (x_1,...,x_d)$. Then just as in Step 1, it is then trivial that there is a constant $M>0$ such that $\chi \le M$ on $\supp \nu$. Then by convexity of $\chi$ we get $\chi \le M$ on $J$ for every $\chi \in \Lambda$, proving the proposition in this simple case. Next if $x \neq a =(a_1,...,a_d)$, recall that in Step 1 we found $\tilde \mu_i, \tilde \nu_i \in \cP(\R)$ such that $\mu_i \le_c \tilde \mu_i \le_c \tilde \nu_i \le_c \nu_i$ and a constant $\theta_i >0$ such that the measure $\eta_i := \theta_i\, \zeta_{a_i \to (y^-_i, y^+_i)}$ satisfies 
$\eta_i \le \tilde \pi_i$ for some $\tilde \pi_i \in {\rm MT}(\tilde \mu_i, \tilde \nu_i)$. Set $\tilde \mu^\otimes := \otimes_i \tilde \mu_i$, $\tilde \nu^\otimes := \otimes_i \tilde \nu_i$ and $ \eta^\otimes := \otimes_i  \eta_i$. Then observe that $\int \chi \,d(\tilde \nu^\otimes-\tilde \mu^\otimes) \le C$ implies $\int \chi \,d(p^y_\# \eta^\otimes -p^x_\# \eta^\otimes) \le C$. Since $\supp (p^x_\# \eta^\otimes) = \{a\}$ and $\supp (p^y_\# \eta^\otimes) = \supp \nu$, we are reduced to the previous case and the local uniform bound follows. 

Now notice we can carry out Step 2 and 3 in this higher dimension case exactly the same way, that is, we can find a countable rectangular martingale measures (say $(\zeta_n)_{n \in \N}$) which covers $I$ and satisfies the appropriate chain condition. We saw in Step 2 that the boundedness property of $\chi$ on each interval $\conv(\supp \zeta_n)$ can propagate along such a  chain, and the argument works the same way in higher dimension, simply observing that the intersecting intervals $\rrbracket y^-_1 + \epsilon , y^+_1 -\epsilon \llbracket \,\, \cap \,\, \rrbracket y^-_2 + \epsilon , y^+_2 -\epsilon  \llbracket \, \neq \emptyset$ in Step 2 are now replaced by intersecting chain of rectangles. Moreover if some marginal measures $\nu_i$ assign positive mass on the boundary of $I_i$, then observe that the rectangles can cover up to such boundaries so that we get the desired boundedness of all $\chi \in \Lambda$ up to the boundary. This completes the proof of the local boundedness.\\

{\bf Step 5.} It remains to prove the convergence property, and it is a direct consequence of the local bound and an application of Arzel{\`a}-Ascoli theorem, as follows: we have shown that every $\chi \in \Lambda$ is uniformly bounded on any compact subset of $J$, so with convexity of $\chi$  and $\chi(a) =0$, $\chi \ge 0$ we deduce that the derivative of every $\chi \in \Lambda$ is also uniformly bounded on any compact subset of $I$. Hence by Arzel{\`a}-Ascoli theorem we deduce that for any sequence $\{ \chi_n \}$ in $\Lambda$ and any compact set $K$ in $I$, there exists a subsequence of $\{ \chi_n \}$ that converges uniformly on $K$. Now by increasing $K$ to $I$ and using diagonal argument, we can find a further subsequence of $\{ \chi_n \}$ that converges pointwise in $I$, and moreover, the convergence is uniform on every compact subset of $I$.

We can deal with the convergence on $J \setminus I$  the same way and let us explain in dimension $2$ to avoid notational difficulty. Suppose, for example,  $J_1 = [a_1, b_1[$, $J_2 = [a_2, b_2[$, and $J = J_1 \times J_2$. First of all, at the corner point $(a_1,a_2)$ the sequence $\{ \chi_n \}$ is bounded, so we can choose a subsequence which converges at $(a_1,a_2)$. Next, on any compact interval in the open line $\{(x,a_2) \,|\, a_1 < x < b_1\}$ the sequence $\{ \chi_n \}$ is also uniformly bounded, so by Arzel{\`a}-Ascoli theorem with diagonal selection again we can find a further subsequence which converges everywhere on $\{(x,a_2) \,|\, a_1 < x < b_1\}$. Finally, we can find a further subsequence which converges on $\{(a_1,y) \,|\, a_2 < y < b_2\}$, and hence the sequence converges everywhere on $J$. This completes the proof of Proposition \ref{convv}.
\end{proof}


\begin{thebibliography}{99}
\bibitem{akp}
L. ~Ambrosio, B. ~Kirchheim, and A. ~Pratelli.
\newblock Existence of optimal transport maps for crystalline norms.
\newblock {\em Duke Mathematical Journal.} Volume 125, Number 2 (2004)  207--241.
 
\bibitem{ap}
L. ~Ambrosio and A. ~Pratelli.
\newblock Existence and stability results in the $L^1$ theory of optimal transportation.
\newblock {\em Optimal transportation and applications. Springer Berlin Heidelberg.},
Volume 1813 of the series Lecture Notes in Mathematics (2003) 123--160.



\bibitem{bch}
M.~Beiglb{\"o}ck, A.M.G.~Cox, and M.~Huesmann.
\newblock Optimal transport and Skorokhod embedding.
\newblock {\em Invent. Math.}, 208 (2017) 327--400.

\bibitem{BeHePe11}
M.~Beiglb{\"o}ck, P.~Henry-Labord{\`e}re, and F.~Penkner.
\newblock Model-independent bounds for option prices: a mass transport approach.
\newblock {\em Finance and Stochastics}, 17 (2013) 477--501.

\bibitem{bj}
M.~Beiglb{\"o}ck and N.~Juillet.
\newblock On a problem of optimal transport under marginal martingale constraints.
\newblock {\em Ann. Probab.}, 44(1) (2016) 42--106.

\bibitem{blo}
M.~Beiglb{\"o}ck, T.~Lim and J.~Ob{\l}{\'o}j.
\newblock Dual attainment for the martingale transport problem.
\newblock {\em Bernoulli}. 25 (2019), no. 3, 1640--1658. 

\bibitem{bnt}
M.~Beiglb{\"o}ck, M. Nutz and N. Touzi.
\newblock  Complete duality for martingale optimal transport on the line.
\newblock {\em Ann. Probab.}, 45(5) (2017) 3038--3074. 

\bibitem{bc}
S.~Bianchini and F.~Cavalletti.
\newblock The Monge Problem for Distance Cost in Geodesic Spaces.
\newblock {\em Comm. Math. Phys.}, Volume 318, Issue 3 (2013) 615--673.

\bibitem{bl78}
D.T. Breeden and R.H. Litzenberger.
\newblock Prices of state-contingent claims implicit in option prices.
\newblock {\em J. Business}, 51(4):621--651 (1978).

\bibitem{br}
Y.~Brenier.
\newblock Polar factorization and monotone rearrangement of vector-valued functions.
\newblock {\em Comm. Pure Appl. Math.}, Volume 44, Issue 4 (1991) 375--417.

\bibitem{cfm}
L.A. Caffarelli, M. Feldman, and R. J. McCann.
\newblock Constructing optimal maps for Monge's transport problem as a limit of strictly convex costs. 
\newblock {\em J. Amer. Math. Soc.}, 15 (2002) 1--26.

\bibitem{ccg16}
G. Carlier, V. Chernozhukov and A. Galichon.
\newblock Vector quantile regression: an optimal transport approach.
\newblock {\em The Annals of Statistics}, 2016, Vol. 44, No. 3, 1165--1192.

\bibitem{cko}
L.~Carraro, N.E.~Karoui and J.~Ob{\l}{\'o}j.
\newblock On Az{\'e}ma-Yor processes, their optimal properties and the Bachelier-drawdown equation.
\newblock   {\em  Ann. Probab.}, Volume 40, Number 1 (2012) 372--400.

\bibitem{cd1}
T.~Champion and L.~De Pascale.
\newblock The Monge problem for strictly convex norms in $\R^d$.
\newblock {\em Journal of the European Mathematical Society.}, 12(6) (2010) 1355--1369.

\bibitem{cd2}
T.~Champion and L.~De Pascale.
\newblock The Monge problem in  $\R^d$.
\newblock {\em Duke Mathematical Journal.}, Volume 157, Number 3 (2011) 551--572.

\bibitem{cfg10}
V. Chernozhukov, I. Fernandez-Val and A. Galichon.
\newblock Quantile and Probability Curves without Crossing.
\newblock  {\em Econometrica} 78(3), pp. 1093--1125 (2010).

\bibitem{cghh17}
V. Chernozhukov, A. Galichon, M. Hallin, and M. Henry.
\newblock Monge-Kantorovich depth, quantiles, ranks and signs.
\newblock {\em The Annals of Statistics}, 2017, Vol. 45, No. 1, 223--256.

\bibitem{cghp21}
V. Chernozhukov, A. Galichon, M. Henry, and B. Pass.
\newblock Identification of hedonic equilibrium and nonseperable simultaneous equations.
\newblock {\em J. Political Econ.}, 
Volume 129, Number 3 (2021).

\bibitem{cms01}
D. Cordero-Erausquin, R. J. McCann, and M. Schmuckenschl{\"a}ger. 
\newblock A Riemannian interpolation inequality {\'a} la Borell,
Brascamp and Lieb.
\newblock {\em Invent. Math.}, 146, (2), (2001), 219--257.

\bibitem{cms06}
D. Cordero-Erausquin, R. J. McCann, and M. Schmuckenschl{\"a}ger. 
\newblock Pr{\'e}kopa-Leindler type inequalities on Riemannian
manifolds, Jacobi fields, and optimal transport.
\newblock {\em Ann. Fac. Sci. Toulouse Math.}, 6, (15), (2006), 613--635.


\bibitem{cot19}
A.M.G.~Cox, J.~Ob{\l}{\'o}j and N.~Touzi.
\newblock The Root solution to the multi-marginal embedding problem: an optimal stopping and time-reversal approach.
\newblock {\em Probab. Theory Related Fields.} 173 (2019), no. 1-2, 211--259.

\bibitem{dlms13}
C. Decker, E.H. Lieb, R.J. McCann and B.K. Stephens.
\newblock  Unique equilibria and substitution effects in a stochastic model of the marriage market.
\newblock {\em J. Econom. Theory} 148 (2013) 778--792.

\bibitem{d18-1}
H. De March.
\newblock Local structure of multi-dimensional martingale optimal transport.
\newblock arXiv preprint. https://arxiv.org/abs/1805.09469

\bibitem{d18}
H. De March.
\newblock Quasi-sure duality for multi-dimensional martingale optimal transport.
\newblock arXiv preprint. https://arxiv.org/abs/1805.01757


\bibitem{dt19}
H. De March and N. Touzi.
\newblock Irreducible convex paving for decomposition of multidimensional martingale transport plans. 
\newblock {\em Ann. Probab.} 47 (2019), No. 3, 1726--1774.

\bibitem{ds1}
Y.~Dolinsky and H.M.~Soner.
\newblock Martingale optimal transport and robust hedging in continuous time.
\newblock {\em Probab. Theory Relat. Fields.}, 160 (2014) 391--427.

\bibitem{eglo}
S.~Eckstein, G.~Guo, T.~Lim, and J.~Ob{\l}{\'o}j.
\newblock Robust Pricing and Hedging of Options on Multiple Assets and Its Numerics.
\newblock {\em SIAM J. Financial Math.}, 12 (2021), no. 1, 158--188.

\bibitem{evans}
L.C.~Evans.
\newblock Partial differential equations and Monge-Kantorovich mass transfer.
\newblock {\em 	Current Developments in Mathematics.}, (1997) International Press.

\bibitem{fkm11}
A. Figalli,Y-H Kim and R.J. McCann.
\newblock When is multidimensional screening a convex program?
\newblock {\em J. Econom. Theory} 146 (2011) 454--478.

\bibitem{GaHeTo11}
A.~Galichon, P.~Henry-Labord{\`e}re, and N.~Touzi.
\newblock A Stochastic Control Approach to No-Arbitrage Bounds Given
  Marginals, with an Application to Lookback Options.
\newblock {\em Annals of Applied Probability.}, Volume 24, Number 1 (2014) 312--336.

\bibitem{gm}
W.~Gangbo and R.J.~McCann.
\newblock The geometry of optimal transportation.
\newblock {\em Acta Math.}, Volume 177, Issue 2 (1996) 113--161. 

\bibitem{GKL2}
N.~Ghoussoub, Y-H. ~Kim, and T. ~Lim.
\newblock Structure of optimal martingale transport plans in general dimensions.
\newblock {\em Ann. Probab.} 47(1): 109--164 (2019). 

\bibitem{GKL3}
N.~Ghoussoub, Y-H. ~Kim, and T. ~Lim.
\newblock Optimal Brownian stopping when the source and target are radially symmetric distributions. 
\newblock {\em SIAM J. Control Optim.} 58 (2020), no. 5, 2765--2789.

\bibitem{gtt1}
G. Guo, X. Tan and N. Touzi.
\newblock On the monotonicity principle of optimal Skorokhod embedding problem.
\newblock {\em SIAM J. Control Optim.}, 54-5 (2016) 2478--2489.

\bibitem{gtt2}
G. Guo, X. Tan and N. Touzi.
\newblock Optimal Skorokhod embedding under finitely-many marginal constraints.
\newblock {\em SIAM J. Control Optim.}, 54-4 (2016) 2174--2201.

\bibitem{Ho98}
D.~Hobson.
\newblock Robust hedging of the lookback option.
\newblock {\em Finance and Stochastics.},  2 (1998) 329--347.

\bibitem{Ho11}
D.~Hobson.
\newblock The {S}korokhod embedding problem and model-independent bounds for
  option prices.
\newblock In {\em Paris-{P}rinceton {L}ectures on {M}athematical {F}inance
  2010}, volume 2003 of {\em Lecture Notes in Math.}, Springer,
  Berlin (2011)  267--318.

\bibitem{HoKl12}
D.~Hobson and M.~Klimmek.
\newblock  Robust price bounds for the forward starting straddle.
\newblock {\em Finance and Stochastics.}, Volume 19, Issue 1 (2014) 189--214.


\bibitem{HoNe11}
D.~Hobson and A.~Neuberger.
\newblock Robust bounds for forward start options.
\newblock {\em Mathematical Finance.}, Volume 22, Issue 1 (2012) 31--56.

\bibitem{K1}
L. V. Kantorovich.
\newblock On the transfer of masses.
\newblock {\em Dokl. Akad. Nauk.}, SSSR 37 (1942) 227--229 (Russian).

\bibitem{K2}
L. V. Kantorovich.
\newblock On a problem of Monge.
\newblock {\em Uspekhi Mat. Nauk.}, 3 (1948) 225--226.

\bibitem{Ke84}
H.~Kellerer.
\newblock Duality theorems for marginal problems.
\newblock {\em Z. Wahrsch. Verw. Gebiete.}, 67(4) (1984) 399--432.


\bibitem{Lim}
T.~Lim.
\newblock Optimal martingale transport between radially symmetric marginals in general dimensions.
\newblock {\em Stochastic Processes and their Applications.} Volume 130, Issue 4 (2020) 1897--1912.


\bibitem{MTW}
X. N. Ma, N. S. Trudinger, and X. J. Wang. 
\newblock Regularity of potential functions of the optimal transportation problem. 
\newblock {\em Arch. Ration. Mech. Anal.}, 177 (2005) no. 2, 151--183.

\bibitem{Mc97}
R. J. McCann.
\newblock A convexity principle for interacting gases.
\newblock {\em Adv. Math.}, 128, 153--179 (1997).

\bibitem{Mo1781}
G.~Monge.
\newblock M\'emoire sur la th\'eorie des d\'eblais et des remblais.
\newblock {\em Histoire de l'acad\'emie {R}oyale des {S}ciences de {P}aris}
  (1781).
  
\bibitem{Obloj}
J.~Ob{\l}{\'o}j.
\newblock  The Skorokhod embedding problem and its offspring.
\newblock {\em  Probability Surveys}, 1 (2004) 321--392.

\bibitem{os17}
J.~Ob{\l}{\'o}j and P. Siorpaes.
\newblock Structure of martingale transports in finite dimensions.
\newblock arXiv preprint. https://arxiv.org/abs/1702.08433

\bibitem{os}
J.~Ob{\l}{\'o}j and P.~Spoida.
\newblock An iterated Az\'{e}ma-Yor type embedding for finitely many marginals.
\newblock {\em Ann. Probab.}, 45(4): 2210--2247 (2017).

\bibitem{p12}
B.~Pass.
\newblock Convexity and multi-dimensional screening for spaces with different dimensions. 
\newblock {\em J. Econom. Theory.} 147 (2012) 2399--2418.

\bibitem{Rock}
R.T.~Rockafellar.
\newblock 
 Characterization of the subdifferentials of convex functions.
 \newblock {\em Pacific J. Math}, 17 (1966) 497--510. 
 
\bibitem{ru}
L.~R{\"u}schendorf. 
\newblock  Construction of multivariate distributions with marginals Given. 
\newblock  {\em Annals of the Institute of Statistical Mathematics.}, 37 (1), (1985) 225--233.

 \bibitem{Sa15}
F. Santambrogio.
\newblock  Optimal transport for applied mathematicians. Calculus of variations, PDEs, and modeling.
\newblock {\em Progress in Nonlinear Differential Equations and their Applications}, 87 Birkhuser/Springer, Cham, (2015)

\bibitem{St65}
V.~Strassen.
\newblock The existence of probability measures with given marginals.
\newblock {\em Ann. Math. Statist.}, 36 (1965) 423-439. 

\bibitem{TW2}
N. S. Trudinger and X. J. Wang. 
\newblock On strict convexity and continuous differentiability of potential functions in optimal transportation. 
\newblock {\em Arch. Ration. Mech. Anal.}, 192 (2009) no. 3, 403--418.

\bibitem{Vi03}
C.~Villani.
\newblock {\em Topics in optimal transportation}, Volume~58, {\em Graduate
  Studies in Mathematics}.
\newblock American Mathematical Society, Providence, RI (2003).

\bibitem{Vi09}
C.~Villani.
\newblock {\em Optimal Transport. Old and New}, Vol. 338, {\em Grundlehren
  der mathematischen Wissenschaften}.
\newblock Springer (2009).

\bibitem{Z}
D.~Zaev.
\newblock On the Monge-Kantorovich problem with additional
linear constraints.
\newblock {\em Mathematical Notes.}, 98(5-6) (2015) 725--741.

\end{thebibliography}
\end{document}